\documentclass{amsart}
\usepackage[all]{xy}

\usepackage{a4wide,amssymb,mathrsfs}

\def\smat#1{\left(\begin{smallmatrix}#1\end{smallmatrix}\right)}

\newtheorem{Pro}{Proposition}[subsection]
\newtheorem{Le}[Pro]{Lemma}
\newtheorem{Th}[Pro]{Theorem}
\newtheorem{Co}[Pro]{Corollary}
\theoremstyle{definition}
\newtheorem{De}[Pro]{Definition}
\theoremstyle{remark}

\def\pa{{\mathscr P}}
\def\ha{{\mathscr H}}

\def\fa{{\mathscr F}}

\let\vt\varTheta
\let\vd\varDelta
\let\th\theta
\let\vf\varphi
\let\sd\ltimes

\let\x\times
\let\ho\simeq

\let\then\Rightarrow
\let\ge\geqslant
\let\le\leqslant

\let\vu\varUpsilon

\let\os\oplus
\let\xto\xrightarrow

\def\0{^{[1]}}
\def\kar{^{\sf Ka}}
\def\mr{\mathring}
\def\wt{\widetilde}

\def\e{_{\mathrm{e}}}

\def\op{^{\mathrm{op}}}
\def\1{^{-1}}

\def\cref#1#2#3{\left(#2\right.\left|\ #3\right)_{#1}}

\def\ka{{\mathscr K}}
\def\ta{{\mathscr T}}
\def\c{{\mathscr C}}
\def\ka{{\mathscr K}}

\def\A{{\mathbb A}}
\def\B{{\mathbb B}}

\def\Z{{\mathbb Z}}

\def\I{{\mathbb I}}
\def\J{{\mathbb J}}

\def\trl{{\sf Triangles}}
\def\cn{{\sf Cand}}
\def\ptr{{{\sf Ptr}}}
\def\trn{{{\sf Triangles}_0}}
\def\Ab{{\sf Ab}}

\def\ker{{\sf Ker}}
\def\Id{{\sf Id}}
\def\id{{\sf id}}
\def\cok{{\sf Coker}}
\def\e{{\sf E}}
\def\HH{{\sf HH}}

\DeclareMathOperator\im{{\mathsf{Im}}}%
\def\hom{{\sf Hom}}%
\DeclareMathOperator\h{{\mathsf{h}}}%

\DeclareMathOperator\cro{{\mathsf{cr}}}%

\numberwithin{equation}{subsection}

\begin{document}

\title{Singular extensions  and triangulated categories}

\author[T. Pirashvili]{Teimuraz  Pirashvili}
\address{
Department of Mathematics\\
University of Leicester\\
University Road\\
Leicester\\
LE1 7RH, UK} \email{tp59-at-le.ac.uk}

\maketitle

\hfill{\emph{Dedicated  in Loving Memory  of Beso Jgarkava}}

\section{Introduction}

In this paper we propose a new look on triangulated categories, which is based on singular
 extensions of additive categories.

Let us recall that if $R$ is a ring and $M$ is a square zero two-sided ideal
of $R$, then $M$ can be considered as a bimodule over the quotient ring
$S=R/M$. Moreover the exact sequence
$$
0\to M\to R\to S\to 0
$$
is a singular extension of the ring $S$ by the bimodule $M$, which
is characterized by an element $e(R)\in\HH^2(S,M)$. Here $\HH^*$
denotes the Hochschild cohomology if $S$ is free as an abelian group
and the Shukla cohomology \cite{sh}, \cite{shukla} in the general
situation. Knowing the triple $(S,M,e(R))$ determines the ring $R$
up to isomorphism. This classical fact admits a straightforward
generalization to preadditive categories known at least from the
work of Mitchell \cite{Mit}.

The above relates to triangulated categories as follows. Let $\ta$
be a triangulated category as it was introduced by Puppe
\cite{puppe}. Thus we do not assume the octahedron axiom of Verdier
\cite{verdier} to hold in $\ta$. We first consider the category
$\ta\0$ of arrows of $\ta$ (see Section \ref{isrebi}). Then for each
morphism $f:A\to B$ of $\ta$ we choose a distinguished triangle:
$$
\xymatrix{A\ar[r]^{f}& B\ar[r]^{u_f} &C_f\ar[r]^{v_f} &A[1]}.
$$
Next we consider the category $\trn(\ta)$ which has the same objects as
$\ta\0$, while morphisms $f\to f'$ in $\trn(\ta)$ are commutative diagrams
$$
\xymatrix{A\ar[r]^{f}\ar[d]^{a}& B\ar[r]^{u_f}\ar[d]^{b}
&C_f\ar[r]^{v_f}\ar[d]^{c} &A[1]\ar[d]^{a[1]}\\
A'\ar[r]^{f'}& B'\ar[r]^{u_{f'}} &C_{f'}\ar[r]^{v_{f'}} &A'[1].}
$$
Consider the functor
$$
\pi:\trn(\ta)\to\ta\0
$$
 which is identity on objects and assigns $(a,b)$ to the triple
$(a,b,c)$. Obviously the functor $\pi$ is identity on objects and surjective
on morphisms. We prove that the kernel of the functor $\pi:\trn({\ta})\to
{\ta}\0$ is a square zero ideal in $\trn(\ta)$ (see Section
\ref{mobservation}). It follows that there exists a bifunctor
$\vt:({\ta}\0)\op\x{\ta}\0\to\Ab$ and a singular extension
$$
0\to \vt\to \trn(\ta)\xto{\pi}\ta\0\to 0.
$$
Hence the category $\trn(\ta)$ and therefore the triangulated
category structure on the category $\ta$ is completely determined by
a bifunctor $\vt$ and the corresponding class
$\vartheta\in\HH^2(\ta\0,\vt)$.

The computation of the bifunctor $\vt$ and of the class
$\vartheta\in \HH^2({\ta}\0,\vt)$ is a hard problem. Of course the
bifunctor $\vt$ and the class $\vartheta$ are not arbitrary and it
is an interesting task to characterize such  pairs
$(\vt,\vartheta)$. In Section \ref{muro} we give a reasonable
solution of this problem. Our first observation is that the
categories involved in our extension possess auto-equivalences
induced by the translation functor of $\ta$. Thus our extension is
in fact a singular $\tau$-extension as it is defined below. Our next
observation is that there exists an easily defined bifunctor $\vd$
(called the Toda bifunctor below), which does not depend on the
triangulated structure at all and is related to the bifunctor $\vt$
via a binatural transformation $\th:\vd\to\vt$ which is an
isomorphism provided one of the arguments is a split morphism of the
category $\ta$. Hence $\vd$ should be considered as a first
approximation of $\vt$. It turns out that in many cases,  but not
always our extension is a pushforward along $\th$. For example this
is  so if $\ta$ is a derived category of a ring (in the classical or
in the brave new algebra sense) and it is not so if $\ta$ is the
triangulated category constructed by Muro \cite{muro}. These facts
lead to the definition of a pseudo-triangulated category in Section
\ref{ppuppe}. We will extend the notion of homology and Massey
triple product from triangulated categories to pseudo-triangulated
categories. Finally in Section \ref{muro} we characterize
triangulated categories among all pseudo-triangulated categories.

\section{Preliminaries}

\subsection{Pre-additive categories}

A category $\A$ together with an abelian group structure on each of the sets
of morphisms $\hom_\A(X,Y)$ is called a \emph{preadditive category} provided
all the composition maps $\hom_\A(Y,Z)\x\hom_\A(X,Y)\to\hom_\A(X,Z)$ are
bilinear maps of abelian groups. Suppose $\A$ and $\B$ are preadditive
categories. A functor $F:\A\to \B$ is said to be an \emph{additive functor}
if the induced maps
$$
\A(X,Y)\to \B(F(X),F(Y)),\ f\mapsto F(f)
$$
are homomorphisms of abelian groups for all objects $X,Y\in\A$.

An \emph{additive category} is a preadditive category $\A$ with zero object
$0$ and such that for all objects $X,Y$ there is given an object $X\oplus Y$
and morphisms
\begin{align*}
i_1:X\to X \oplus Y,\ \ \ &i_2:Y\to X\oplus Y,\\
r_1:X\oplus Y \to X,\ \ \ &r_2:X\oplus Y\to Y
\end{align*}
with $r_1i_1=\id_X$, $r_2i_2=\id_Y$, $r_1i_2=0$, $r_2i_1=0$ and
$i_1r_1+i_2r_2=\id_{X\oplus Y}$. The object $X \oplus Y$ is called
\emph{direct sum} of $X$ and $Y$ in $\A$. It follows that $X\oplus Y$
together with $i_1$ and $i_2$ is a coproduct of $X$ and $Y$ and $X\oplus Y$
together with $r_1$ and $r_2$ is a product of $X$ and $Y$. The following fact
is well known.

\begin{Le}
For additive categories $\A$ and $\B$, a functor $F:\A\to\B$ is additive iff
for all objects $X_1$, $X_2$ of the category $\A$ the canonical map
$$
(F(r_1),F(r_2)):F(X_1\os X_2)\to F(X_1)\os F(X_2)
$$
is an isomorphism.
\end{Le}

\

\subsection{Split idempotent and split morphisms}
Let $e:A\to A$ be an endomorphism. If $e^2=e$ then $e$ is called
\emph{idempotent}. If $e$ is an idempotent in an additive category
$\A$ then $\id_A-e$ is also an idempotent. For any objects $X_1$ and
$X_2$ of an additive category $\A$, the morphism $e=i_1r_1:X_1\os
X_2\to  X_1\os X_2$ is an idempotent. An idempotent $e:A\to A$ is
called \emph{split} if there are arrows (called \emph{splitting
data}) $a:A\to B$ and $b:B\to A$, such that $e=ba$ and $ab=\id_B$.
An additive category $\A$ is called \emph{Karoubian} provided all
idempotents split, which is the same to require as that all
idempotents have kernels (or cokernels).

 A morphism $p:X\to Y$ of an
additive category is called a \emph{splittable epimorphism} if there
exists a morphism $j:Y\to X$ such that $pj=\id_Y$. For example the
canonical projection $r:A\os B\to B$ is splittable. Morphisms
isomorphic to such projections are called \emph{split epimorphisms}.
If $\A$ is Karoubian then any splittable epimorphism is actually a
split epimorphism.

Dually a morphism $i:X\to Y$ is  called a \emph{splittable monomorphism} if
there exists a morphism $r:Y\to X$ such that $ri=\id_X$. For example the
canonical inclusion $i:A\to A\os B$ is splittable. Morphisms isomorphic to
such inclusions are called \emph{split monomorphisms}. If $\A$ is Karoubian
then any splittable monomorphism is actually a split monomorphism.

More generally a morphism $f:X\to Y$ is called \emph{splittable} if there
exist a morphism $s:Y\to X$ such that $fsf=f$. Examples of splittable
morphisms are splittable epimorphisms, splittable monomorphisms and
idempotents. Morphisms of the form $\smat{1&0\\0&0}:X\oplus X'\to X\oplus
X''$ are split morphisms. Morphisms isomorphic to such a morphism are called
\emph{split morphisms}. If $\A$ is a Karubian category then any splittable
morphism $f$ is actually a split morphism, i.~e. it can be represented as a
composite $ir$, where $r$ is a split epimorphism and $i$ is a split
monomorphism.

\subsection{Subfunctors of additive functors and the second cross-effect}

Let $\A$ be an additive category. Let $F:\A\to\Ab$ be a functor with
$F(0)=0$. The \emph{second cross-effect} functor of $F$ is a bifunctor
$\cro_2(F):\A\x\A\to\Ab$ defined by
$$
\cro_2(F)(X_1,X_2):=\ker((F(p_1),F(p_2)):F(X_1\os X_2)\to F(X_1)\os F(X_2)).
$$
Thus a functor $F$ is additive iff $F(0)=0$ and $\cro_2(F)=0$.

The proof of the following fact is an easy exercise on diagram chase and is
left to the reader.

\begin{Le}\label{qvesenijirime}
For any short exact sequence of functors
$$
0\to F_1\to F\to F_2\to0
$$
one has a short exact sequence of bifunctors:
$$
0\to\cro_2(F_1)\to\cro_2(F)\to\cro_2(F_2)\to0
$$
In particular any subfunctor of an additive functor is also additive.
\end{Le}\qed

\subsection{Ideals and quotient categories}

An \emph{ideal} $\I$ of $\A$ is a subbifunctor of the bifunctor
$$
\hom_{\A}(-,-):\A\op\x\A\to\Ab.
$$
It follows from Lemma \ref{qvesenijirime} that $\I$ is biadditive. If $\A$
and $\B$ are additive categories and $F:\A\to \B$ is an additive functor, one
denotes by $\ker(F)$ the ideal of $\A$ consisting of morphisms $f:A\to B$
such that $F(f)$ is a zero morphism in $\B$.

If $\I$ is an ideal of $\A$, then one can form the quotient category $\A/\I$,
which has the same objects as $\A$, while morphisms in $\A/\I$ are given by
$$
\hom_{\A/\I}(A,B):=\hom_{\A}(A,B)/\I(A,B).
$$
One has the canonical additive functor $Q:\A\to \A/\I$. It is clear that
$\ker(Q)=\I$. Any additive functor $F:\A\to\B$ factors through the category
$\A/\ker(F)$.

\subsection{Nilpotent and square zero ideals}

Let $\I$ and $\J$ be ideals of $\A$. For all object $A$ and $B$ we let
$\I\J(A,B)$ be the set of all products $fg$, where $f\in \I(C,B)$ and $g\in
\J(A,C)$, for some $C$. We claim that $\I\J(A,B)$ is a subgroup of $\A(A,B)$.
Indeed, if $f\in\I(C,B)$, $g\in\J(A,C)$ and $f'\in\I(C',B)$, $g'\in\J(A,C')$,
then $fg+f'g'=f''g''$, where $f''=(f,f'):C\os C'\to B$ and $g''=\smat{g\\
g'}: A\to C\os C'$. We have $\J(A,C\os C')=\J(A,C)\os\J(A,C')$ by Lemma
\ref{qvesenijirime}. Since $g\in\J(A,C)$ and $g'\in\J(A,C')$, it follows that
$g''\in \J(A,C\os C')$. Similarly $f''\in\I(C\os C',B)$, hence the claim. It
is clear that $\I\J$ is a subbifunctor of $\J$ and $\J$. Hence it is an
ideal.

Having defined product of ideals, one can talk about powers $\I^n$ of an
ideal $\I$. An ideal $\I$ is \emph{nilpotent} if $\I^n=0$ for some $n$. Of
special interest are ideals with $\I^2=0$. They are called \emph{square zero}
ideals. We have the following easy but useful fact.

\begin{Le}\label{i2=0}
For a square zero ideal $\I$ of $\A$ the bifunctor $\I:\A\op\x\A\to\Ab$
factors through the quotient category $\A/\I$ in an unique way.
\end{Le}\qed

This result can be used to prove the following simple result.

\begin{Le}\label{refisonil}
Let $\I$ be a nilpotent ideal of an additive category $\A$. Then the quotient
functor $Q:\A\to\A/\I$ reflects isomorphisms and yields an isomorphism of
monoids of isomorphism classes ${\rm Iso}(\A)\cong{\rm Iso}(\A/\I)$.
\end{Le}

\begin{proof}
The last statement follows from the previous one, because $Q$ is identity on
objects and surjective on morphisms. To prove the first statement, it
suffices to assume that $\I^2=0$. Let $f:A\to B$ be a morphism, such that
$Q(f)$ is an isomorphism. Thus there exists $g:B\to A$ such that
$a=gf-\id_A\in \I(A,A)$ and $b=fg-\id_B\in \I(B,B)$. Since $\I$ as a
bifunctor factors through $Q$, it follows that the map $\I(B,A)\to\I(A,A)$
given by $x\mapsto xf$ is an isomorphism. Thus there exists a $c\in\I(B,A)$
with $cf=a$. Now we put $g_1=g-c$. Then $g_1f=gf-cf=\id_A$, which shows that
$f$ is a splittable monomorphism. A similar argument shows that $f$ is a
splittable epimorphism, hence an isomorphism. Thus $Q$ reflects isomorphisms
and we are done.
\end{proof}

\subsection{Singular extensions of additive categories}

Let $\B$ be an additive category and let $D:\B\op\x \B\to \Ab$ be a
bifunctor. A \emph{singular extension}
$$
0\to D\xto i\A\xto F\B\to 0
$$
of $\B$ by the  bifunctor $D$ is the following data:

\begin{enumerate}
\item
An additive category $\A$ and an additive functor $F:\A\to\B$, such that
$\ker(F)$ is a square zero ideal and the canonical functor $\A/\ker(F)\to\B$
is an isomorphism of categories;
\item
an isomorphism of bifunctors $i:D(F(\cdot),F(\cdot))\to\ker(F)$.
\end{enumerate}

\subsection{Semidirect product}

Let $\B$ be an additive category and let $D:\B\op\x\B\to\Ab$ be a bifunctor.
The \emph{semidirect product} (compare with \cite{BW}) of $\B$ by $D$ is the
category $\B\sd D$ which has the same objects as $\B$. Morphisms $A\to B$ in
$\B\sd D$ are pairs $(f,a)$, where $f:A\to B$ is a morphism in $\B$ and $a\in
D(A,B)$. Composition is defined by
$$
(f,a)\circ (g,b)=(fg,f_*(b)+g^*(a))
$$
Let $\I$ be the class of all morphisms of the form $(0,a)$. Then $\I^2=0$,
$\A/\I\cong\B$, where $\A=\B\sd D$ and $i:D\to\I$ is an isomorphism of
bifunctors, given by $i(a)=(0,a)$. Conversely, if
$$
0\to D\xto i\A\xto F\B\to0
$$
is a singular extension and $F$ has a section, then $\A\cong\B\sd D$.

\subsection{Cohomology and singular extensions}

The reader familiar with the Hoch\-schild cohomology and especially with
relations between the second Hochschild cohomology and singular extensions of
rings might wonder whether there is a cohomology theory which in dimension
two would classify singular extensions of a small additive category $\B$ by a
bifunctor $D:\B\op\x\B\to\Ab$. In fact such cohomology does exist and it is
an obvious extension of the Shukla cohomology of rings \cite{sh},
\cite{shukla} to small preadditive categories.

As a matter of fact, let us mention here that there exists also
Baues-Wirsching cohomology \cite{BW} which is defined for all small (maybe
non-preadditive) categories. For additive categories the second Shukla
cohomology and the second Baues-Wirsching cohomology $H^2(\A,D)$ are
isomorphic. This follows from \cite{BW}, together with Proposition 3.4 of
\cite{JP}, which shows that any linear extension \cite{BW} of an additive
category by an additive bifunctor is again an additive category. It must be
mentioned that even for additive categories Shukla and Baues-Wirsching
cohomologies are not isomorphic in dimensions $\ge3$.

\subsection{Puppe triangulated categories}

Let $\ta$ be an additive category with an autoequivalence $A\mapsto
A[1]$. A \emph{candidate triangle} in $\ta$ is a diagram $$ X\to
Y\to Z\to X[1].$$ A morphism from  a candidate triangle $ X\to Y\to
Z\to X[1]$ to a candidate triangle $ X'\to Y'\to Z'\to X'[1]$ is a
commutative diagram in $\ta$:
$$
\xymatrix{
X\ar[r]\ar[d]^a&Y\ar[r]\ar[d]^b&Z\ar[r]\ar[d]^c&X[1]\ar[d]^{a[1]}\\
X'\ar[r]&Y'\ar[r]&Z'\ar[r]&X'[1] }
$$
We let $\cn$ be the category of candidate triangles. A candidate
triangle  $X\to Y\to Z\to X[1]$ is \emph{acyclic} provided the
sequence of abelian groups
$$\cdots \to \hom_\ta(X[1],A)\to \hom_\ta(Z,A)\to \hom_\ta(Y,A)\to
\hom_\ta(X,A)\to \cdots$$ is exact for any object $A\in \ta$.

 A \emph{Puppe triangulated category}
structure, or simply triangulated category structure on $\ta$ is
given by a collection of diagrams, called \emph{distinguished
triangles}, of the form
$$
X\to Y\to Z\to X[1]
$$
such that
\begin{itemize}
\item[TR1)]
Any candidate triangle isomorphic to a distinguished triangle in
$\cn$ is a distinguished triangle.
\item[TR2)]
Any diagram of the following form is a distinguished triangle:
$$
X\xto{\id_X}X\to0\to X[1]
$$
\item[TR3)]
If
$$
X\xto fY\xto gZ\xto hX[1]
$$
is a distinguished triangle, then
$$
Y\xto gZ\xto hX[1]\xto{-f[1]}Y[1]
$$
is also a distinguished triangle.
\item[TR4)]
For any morphism $f:X\to Y$ there is a distinguished triangle of the form
$$
X\xto fY\to Z\to X[1].
$$
\item[TR5)]
Suppose we have a diagram
$$
\xymatrix{
X\ar[r]\ar[d]^a&Y\ar[r]\ar[d]&Z\ar[r]&X[1]\ar[d]^{a[1]}\\
X'\ar[r]&Y'\ar[r]&Z'\ar[r]&X'[1]
}
$$
in which the rows are distinguished triangles and the left rectangle
commutes. Then there exists a morphism $Z\to Z'$ making the diagram
$$
\xymatrix{
X\ar[r]\ar[d]^a&Y\ar[r]\ar[d]&Z\ar[r]\ar[d]&X[1]\ar[d]^{a[1]}\\
X'\ar[r]&Y'\ar[r]&Z'\ar[r]&X'[1]
}
$$
commute.
\end{itemize}

A category equipped with  a triangulated structure is called a
\emph{triangulated category}. We let $\trl(\ta)$ be the full
subcategory of $\cn$ formed by distinguished triangles.

Let $\ta$ be a triangulated category. An additive functor $\h:\ta\to
\Ab$ is called \emph{homology} if, whenever
$$
X\xto{f} Y\to Z\to X[1]
$$
is a distinguished triangle, the sequence
$$
\h(X)\to \h(Y)\to\h(Z)\to\h(X[1])
$$
is exact. Then the sequence
$$
\cdots\to\h^n(X)\to\h^n(Y)\to\h^n(Z)\to\h^{n+1}(X)\to\cdots
$$
is also exact, where $\h^n(X)=\h(X[n])$.

It is well known that the functors $\hom_\ta(X,-)$ and $\hom_\ta(-,X)$ are
homologies. In particular if $X\xto fY\xto gZ\to X[1]$ is a distinguished
triangle and $h:Y\to V$ is a morphism such that $hf=0$, then $h$ factors
through $g$.

\section{Singular extensions and triangulated
categories}

\subsection{Category of arrows}\label{isrebi}

Let $[1]$ be the category associated to the ordered set $0<1$. For
any category $\c$ we let $\c\0$ be the category of functors $[1]\to
\c$. Thus $\c\0$ is the category of arrows of $\c$. For a morphism
$f:A\to B$ of the category $\c$ considered as an object of the
category $\c\0$ we use the notation $\mr f$ and the word ''arrow''
to denote the same morphism considered as an object of the category
$\c\0$. Hence objects of $\c\0$ are arrows $\mr f$, where $f:A\to B$
is a morphisms of $\c$, while morphisms $\mr f\to\mr f'$ are pairs
of morphisms $(a:A\to A', b:B\to B')$ in $\c$ such that the diagram
$$
\xymatrix{
A\ar[r]^{f}\ar[d]^{a}& B\ar[d]^{b} \\
A'\ar[r]^{f'}& B' }
$$
commutes.

For any object $A$ of $\A$ we write $\id_A$ for the identity morphism in $\A$
and use $\Id_A$ for the corresponding arrow considered as an object of
$\c\0$. Hence $\Id_A=\mr\id_A$. Assume now that $\c$ has a zero object. In
this case we use the following notations. For an object $A$ in $\c$ we denote
by $^A!$ (resp. $!_A$) the object of $\A\0$ corresponding to the unique
morphism $0\to A$ (resp. $A\to0$) in $\c$.

The functors
$$
\Id_{?},\, !_{?},\, ^{?}!:\c\to \c\0
$$
are full embeddings.

\subsection{The main observation}\label{mobservation}

Let $\ta$ be a triangulated category. For each morphism $f:A\to B$ of $\ta$
we choose a distinguished triangle
\begin{equation}\label{samk}
\xymatrix{
A\ar[r]^{f}& B\ar[r]^{u_f} &C_f\ar[r]^{v_f} &A[1]},
\end{equation}
where $A\mapsto A[1]$ is the translation functor. One of the axioms of
triangulated categories asserts that such choice is always possible. Now we
consider the category $\trn(\ta)$, whose objects are morphisms of $\ta$, thus
the same as of the category $\ta\0$. For a morphism $f:A\to B$ we let $[f]$
be the corresponding object of the category $\trn(\ta)$. The morphisms
$[f]\to[f']$ in the category ${\trn}({\ta})$ are triples of morphisms
$(a:A\to A',b:B\to B',c:C_f\to C_{f'})$ of the category $\ta$ such that the
diagram
$$
\xymatrix{
A\ar[r]^{f}\ar[d]^{a}& B\ar[r]^{u_f}\ar[d]^{b} &C_f\ar[r]^{v_f}\ar[d]^{c} &A[1]\ar[d]^{a[1]}\\
A'\ar[r]^{f'}& B'\ar[r]^{u_{f'}} &C_{f'}\ar[r]^{v_{f'}} &A'[1]
}
$$
is commutative. Thus we have full subcategories $\trn(\ta)\subset
\trl(\ta)\subset \cn$. It is clear that the first inclusion
$\trn(\ta)\subset \trl(\ta)$ is an equivalence of categories.
Moreover the category $\trl(\ta)$ can be reconstructed from
$\trn(\ta)$ as follows: A candidate triangle belongs to $\trl(\ta)$
iff if it is isomorphic (in $\cn$) to an object of $\trn(\ta)$.

We let
$$
\pi:\trn(\ta)\to\ta\0
$$
be the functor which is the identity on objects (thus $\pi([f])=\mr{f}$) and
assigns $(a,b)$ to the triple $(a,b,c)$. Another axiom of triangulated
categories asserts that the functor $\pi$ is surjective on morphisms.

\begin{Le}\label{conrep}
For arbitrary object $X$ in a triangulated category $\ta$ and arbitrary
morphism $f:A\to B$, there exist isomorphisms
$$
\hom_\ta(C_f,X)\cong \hom_{\trn(\ta)}([f],!_X)
$$
and
$$
\hom_\ta(X,C_f[-1])\cong \hom_{\trn(\ta)}(^X!,[f]).
$$
These isomorphisms are natural in $X\in\ta$ and in $f\in\trn(\ta)$.
\end{Le}

\begin{proof}
We prove the first isomorphism, second being similar. A morphism
$$
\xymatrix{
A\ar[r]^{f}\ar[d]^{0}& B\ar[r]^{u_f}\ar[d]^{b} &C_f\ar[r]^{v_f}\ar[d]^{c} &A[1]\ar[d]^{0}\\
0\ar[r]^{0}& X\ar[r]^{\id} &X\ar[r]^{v_{0}}&0
}
$$
is uniquely determined by $c$, which might be arbitrary. This implies the
result.
\end{proof}

The following easy but extremely important fact is new.

\begin{Le}\label{hemokargo}
The kernel of the functor $\pi:\trn(\ta)\to\ta\0$ is a square zero ideal.
\end{Le}

\begin{proof}
Consider the following commutative diagram in $\ta$
$$
\xymatrix{
A\ar[r]^{f}\ar[d]^{0}& B\ar[r]^{u_f}\ar[d]^{0}&C_f\ar[r]^{v_f}\ar[d]^{c} &A[1]\ar[d]^{0}\\
A'\ar[r]^{f'}\ar[d]^0& B'\ar[r]^{u_{f'}}\ar[d]^0 &C_{f'}\ar[r]^{v_{f'}}\ar[d]^{c'}&A'[1]\ar[d]^0\\
A''\ar[r]^{f''}& B''\ar[r]^{u_{f''}} &C_{f''}\ar[r]^{v_{f''}} &A''[1]
}
$$
where rows are  distinguished triangles. We have to prove that $c'c=0$. Since $c'u_{f'}=0$,
there exist a morphism $d':A'[1]\to C_{f''}$ such that $c'=d'v_{f'}$. Hence $$c'c=d'v_{f'}c=d'0v_f=0.$$
\end{proof}

\begin{Co}\label{333sed}
There exists a well-defined bifunctor $\vt_\ta$
\begin{equation}\label{hemiklasi}
\vt_\ta:(\ta\0)\op\x\ta\0\to\Ab
\end{equation}
such that
$$
\vt_\ta(\mr f,\mr{f'})=\{c:C_f\to C_{f'}\mid cu_f=0,v_{f'}c=0\}.
$$
The category $\trn(\ta)$ is a singular extension of the category $\ta\0$ by
the bifunctor $\vt_\ta$,
\begin{equation}\label{hemigafa}
0\to\vt\to\trn(\ta)\xto\pi\ta\0\to0.
\end{equation}
The class $\vartheta$ of the singular extension (\ref{hemiklasi}) in
$\HH^2(\ta,\vt)$ is independent of the choices of distinguished triangles
(\ref{samk}). Hence the triangulated category structure on the category $\ta$
is completely determined by the bifunctor $\vt$ and the class $\vt_\ta$.
\end{Co}

\subsection{Categories with translation}

Let $\A$ be an additive category. A \emph{translation} on a category $\A$ is
an autoequivalence $\A\to \A$; if such a translation is fixed, then we say
that $\A$ is a \emph{category with translation} or \emph{$\tau$-category}. An
evaluation of the translation functor on an object $A$ is denoted by $A[1]$
and is called translation of $A$. Moreover, for any object $A$ we choose an
object $A[-1]$ together with an isomorphism $(A[-1])[1]\cong A$. Then
$A\mapsto A[-1]$ can be extended as a functor $(\cdot)[-1]:\A\to\A$ in a
unique way. If $n$ is an integer, then one has objects $A[n]$ defined by
induction: $A[n+1]=(A[n])[1]$ if $n\ge1$, $A[0]=A$ and $A[n-1]=(A[n])[-1]$ if
$n\le-1$. Sometimes  we write $\tau(A)$ instead of $A[1]$. Of course in this
case we write $\tau^n(A)$ instead of $A[n]$ as well.

Let $\A$ and $\B$ be categories with translation. A \emph{translation
preserving functor}, or \emph{$\tau$-functor} is an additive functor $F:\A\to
\B$ such that $F(A[1])=(F(A))[1]$ for all $A$.

Let $\I$ be an ideal in a $\tau$-category $\A$. We will say $\I$ is a
\emph{$\tau$-ideal} if for all objects $A$ and $B$ the isomorphism
$\A(A,B)\to\A(A[1],B[1])$, $f\mapsto f[1]$, restricts to an isomorphism
$\I(A,B)\to\I(A[1],B[1])$. In this case the quotient category $\A/\I$ carries
a $\tau$-category structure and the quotient functor $\A\to\A/\I$ is a
$\tau$-functor. Conversely, if $F:\A\to\B$ is a $\tau$-functor, then
$\ker(F)$ is a $\tau$-ideal.

\subsection{Koszul translation}

For a morphism $f:X\to Y$ in a $\tau$-category $\A$ one puts:
$$
\tau(\mr f)=(-f[1]:X[1]\to Y[1]).
$$
Moreover, if  $f':X'\to Y'$ is another morphism of the category $\A$ and
$(x:X\to X',y:Y\to Y')$ is a morphism $\mr f\to\mr{f'}$ in the category
$\A\0$, then one puts
$$
\tau(x,y)=((x[1],[y]):\tau(\mr{f})\to \tau(\mr{f'})).
$$
In this way one gets a translation $\tau:\A\0\to\A\0$ called the \emph{Koszul
translation}.

Let $\ta$ be a triangulated category. Then $\trn(\ta)$ also possesses a
Koszul translation, which on objects is given by the same rule
$$
\tau([f])=(-f[1]:X[1]\to Y[1]),
$$
while on morphisms it is given by $\tau(x,y,z)=(x[1],y[1],z[1])$. Here
$(x[1],y[1],z[1])$ is the following morphism in $\trn(\ta)$:
$$
\xymatrix{%
X[1]\ar[r]^{-f[1]}\ar[d]^{x[1]}&Y[1]\ar[r]^{-u_f[1]}\ar[d]^{y[1]}&C_f[1]\ar[r]^{-v_f[1]}
\ar[d]^{z[1]} &X[2]\ar[d]^{x[2]}\\
X'[1]\ar[r]^{-f'[1]}& Y'[1]\ar[r]^{-u_{f'}[1]}&C_{f'}[1]\ar[r]^{-v_{f'}[1]}&A'[2]%
}
$$
It is clear that $\pi:\trn(\ta)\to\ta\0$ is a $\tau$-functor.

\subsection{$\tau$-bifunctors and singular $\tau$-extensions}

Let $\A$ be a $\tau$-category. A \emph{$\tau$-bifunctor} on $\A$ is a
bifunctor $D:\A\op\x\A\to\Ab$ together with a system of isomorphisms
$$
t_{A,B}:D(A,B)\to D(A[1],B[1]),\ A,B\in\A,
$$
which are natural in $A$ and $B$. If $D$ and $D'$ are two $\tau$-bifunctors,
then a natural transformation $\xi:D\to D'$ of bifunctors is called a
\emph{$\tau$-transformation} provided the following diagram commutes:
$$
\xymatrix{
D(A,B)\ar[d]_{\xi(A,B)}\ar[r]^{t_{A,B}}&D(A[1],B[1])\ar[d]^{\xi(A[1],B[1])}\\
D'(A,B)\ar[r]^{t'_{A,B}}&D'(A[1],B[1])
}
$$

For example the bifunctor $\hom_\A(\cdot,\cdot)$ is a $\tau$-bifunctor, where
$t_{A,B}(f)=f[1]$. Moreover, if $\I$ is a $\tau$-ideal, then it is a
$\tau$-subbifunctor of $\hom_\A$.

A singular extension
$$
0\to D\xto i\B\xto p\A\to0
$$
of a $\tau$-category $\A$ by a $\tau$-bifunctor $D$ is called a
\emph{singular $\tau$-extension} if $p$ is a $\tau$-functor and $i$ yields an
isomorphism $D\to\ker(p)$ of $\tau$-bifunctors over $\A$.

One easily sees that the singular extension \ref{hemigafa} is in fact a
singular $\tau$-extension, where $\ta\0$ and $\trn(\ta)$ are equipped with
Koszul translations. Here a $\tau$-bifunctor structure on $\vt$, i.~e.
isomorphisms
$$
t_{\mr f,\mr{f'}}:\vt(\mr f,\mr{f'})\to\vt(-\mr f[1],-\mr{f'}[1])
$$
are induced by $c\mapsto c[1]$, for any $c:C_f\to C_{f'}$ with $cu_f=0=v_{f'}c$.

\subsection{Toda bifunctor}

Let $\A$ be a category with translation. For morphisms $f:A\to B$ and
$f':A'\to B'$ we consider the homomorphism of abelian groups
$$
\phi_{f,f'}:\hom_\A(A[1],A')\oplus\hom_\A(B[1],B')\to\hom_\A(A[1],B')
$$
given by
$$
\phi_{f,f'}(g,h)=f'_*(g)-(f[1])^*(h)=f'\circ g-h\circ(f[1]).
$$
Here $g:A[1]\to A'$ and $h:B[1]\to B'$ are morphisms of $\A$.

The \emph{Toda bifunctor} $\vd_\A$, or simply $\vd$ is a bifunctor
$$
\vd:(\A\0)\op\x\A\0\to\Ab
$$
given by
$$
\vd(\mr f,\mr{f'}):=\cok(\phi_{f,f'})=\frac{\hom_{\A}(A[1],B')}
{f'_*\hom_{\A}(A[1],A')-f^*\hom_{\A}(B[1],B')},
$$
where $f:A\to B$ and $f':A'\to B'$ are morphisms in $\A$.

The following lemma is straightforward.

\begin{Le}\label{d=0}
Let $\A$ be a category with translation. For any object $X\in\c$ and any
morphism $f:A\to B$ one has
\begin{align*}
\vd(\Id_X,\mr f)&=0,\\
\vd(\mr f,\Id_X)&=0,\\
\vd(!_X,\mr f)&=0,\\
\vd(\mr f,^X!)&=0,\\
\vd(^X!,\mr f)&=\cok(\hom_\A(X[1],A)\xto{f_*}\hom_\A(X[1],B)),\\
\vd(\mr f,!_X)&=\cok(\hom_\A(B[1],X)\xto{-f^*}\hom_\A(A[1],X)).
\end{align*}
\end{Le}

There is a $\tau$-bifunctor structure on $\vd$ which we will use throughout.
The isomorphisms
$$
t_{\mr f,\mr{f'}}:\vd(\mr f,\mr{f'})\to\vd(-\mr f[1],-\mr{f'}[1])
$$
are induced by $a\mapsto a[1]$.

\subsection{Natural transformation $\theta$}

Let $\ta$ be a triangulated category. Then we have two $\tau$-bifunctors
$$
\vt_\ta,\vd_\ta:(\ta\0)\op\x\ta\0\to\Ab.
$$
It must be noticed that the Toda bifunctor depends only on the translation
structure, while the bifunctor $\vt$ depends on the choice of the class of
distinguished triangles. We now define the $\tau$-transformation
$$
\th_\ta:\vd_\ta\to\vt_\ta
$$
as follows.

Let $f:A\to B$ and $f':A'\to B'$ be morphisms in $\ta$. For any morphism
$x:A[1]\to B'$ we have the following morphism of distinguished triangles:
$$
\xymatrix{
A\ar[r]^f\ar[d]^0&B\ar[r]^{u_f}\ar[d]^0&C_f\ar[r]^{v_f}\ar[d]^{c_x}&A[1]\ar[d]^0\\
A'\ar[r]^{f'}&B'\ar[r]^{u_{f'}}&C_{f'}\ar[r]^{v_{f'}}&A'[1],
}
$$
where $c_x=u_{f'}xv_f$. One  easily sees that the assignment $x\mapsto
(0,0,c_x)$ yields the homomorphism
$\th(\mr{f},\mr{f'}):\vd(\mr{f},\mr{f'})\to \vt(\mr{f},\mr{f'})$, hence a
natural transformation $\th:\vd\to\vt$.

\begin{Pro}
Let $f:A\to B$ and $f':A'\to B'$ be morphisms of $\ta$. If $f$ or $f'$ is
splittable, then $\th(\mr f,\mr{f'}):\vd(\mr f,\mr{f'})\to \vt(\mr
f,\mr{f'})$ is an isomorphism.
\end{Pro}

\begin{proof}
It is well known that for any triangulated category $\ta$ the
Karoubian completion $\ta\kar$ has also a triangulated category
structure and the inclusion functor $i$ is a morphism of
triangulated categories. Hence without loss of generality we may
assume that all idempotents split in $\ta$. By duality it suffices
to consider the case, when $f$ is splittable. Since $\ta\0$ is an
additive category, $\th$ is a transformation of additive bifunctors
and any splittable morphism considered as an object of $\ta\0$ is
isomorphic to a direct sum of objects of the form $\Id_A$, $!_B$ or
$^C!$, we have to consider three cases $\mr f=\Id_X$, $\mr f=!_X$
and $\mr f=^X!$. In the first case we have $C_f=0$ and therefore
both groups $\vd(\mr f,\mr{f'})$ and $\vt(\mr f,\mr{f'})$ are
trivial. If $\mr f=!_X$, then we have already shown that
$\vd(\mr{f},\mr{f'})=0$. On the other hands if
$$
\xymatrix{
0\ar[r]^{0}\ar[d]^0&X\ar[r]^\id\ar[d]^0&X\ar[r]^0\ar[d]^c&0\ar[d]^0\\
A'\ar[r]^{f'}&B'\ar[r]^{u_{f'}}&C_{f'}\ar[r]^{v_{f'}}&A'[1]
}
$$
is a morphism in $\trn(\ta)$ then $c=0$. Hence $\vt(\mr f,\mr{f'})=0$ as
well. Now consider the case, when $\mr f=^X!$. Let $c:X[1]\to C_{f'}$ be a
morphism in $\ta$, then
$$
\xymatrix{
X\ar[r]^0\ar[d]^0&0\ar[r]^0\ar[d]^0&X[1]\ar[r]^{-\id}\ar[d]^c&X[1]\ar[d]^0\\
A'\ar[r]^{f'}&B'\ar[r]^{u_{f'}}&C_{f'}\ar[r]^{v_{f'}}&A'[1]
}
$$
is a morphism of distinguished triangles iff
$$
c\in\ker\left(\hom_\A(X[1],C_{f'})\to\hom_\A(X[1],A'[1])\right).
$$
But the last group is isomorphic to
$\cok\left(\hom_\A(X[1],A')\to\hom_\A(X[1],B')\right)=\vd(\mr
f,\mr{f'})$ and we are done.
\end{proof}

\section{Pseudo-triangulated categories}

\subsection{Definition and examples}\label{ppuppe}

Let $\pa$ be an additive category with translation $A\mapsto A[1]$.

\begin{De}
We will say that there is given a \emph{pseudo-triangulated}
category structure on $\pa$ if there is given a singular
$\tau$-extension
$$
0\to\vu\xto i\ptr\xto p\pa\0\to0
$$
of $\pa\0$ by a $\tau$-bifunctor
$$
\vu:(\pa\0)\op\x\pa\0\to\Ab
$$
together with a $\tau$-transformation $\vf:\vd\to \vu$ from the Toda
bifunctor to $\vu$ such that $\vf(\mr f,\mr{f'}):\vd(\mr
f,\mr{f'})\to \vu(\mr f,\mr{f'})$ is an isomorphism provided $f$ or
$f'$ is splittable. If additionally $\vf$ is isomorphic then we say
that $\pa$ is equipped with a \emph{Toda pseudo-triangulated
category} structure.
\end{De}

Abusing notation we will say that $\pa$ is a pseudo-triangulated
category provided such a structure is given.

We have already seen that if $\pa=\ta$ is a triangulated category then the
extension
$$
0\to\vt\to\trn(\ta)\xto\pi\ta\0\to0
$$
together with the transformation $\th:\vd\to\vt$  gives rise to a
pseudo-triangulated structure on $\ta$. We refer to this example as
the \emph{pseudo-triangulated category associated} to a triangulated
category $\ta$.

Unlike the triangulated category structure, any $\tau$-category can
be equipped with the structure of a pseudo-triangulated category:
one can take $\vu=\vd$, and define $\ptr$ to be the semidirect
product of $\pa\0$ with $\vd$, or one can take any other singular
$\tau$-extension of $\pa\0$ by $\vd$. Thus triangulated category
structures on a given category might be really different.

Let $\pa$ be a pseudo-triangulated category. Objects of the category
$\ptr$ are the same as of $\pa\0$, i.~e. they are still arrows, but
now called \emph{pseudo-triangles}. If a morphism $f$ is considered
as a pseudo-triangle, we use the notation $[f]$ instead of $f$.
Assume $[f]$ and $[f']$ are  pseudo-triangles. Then one has the
exact sequence of abelian groups
\begin{equation}\label{des}
0\to\vu(\mr f,\mr{f'})\xto i\hom_\ptr([f],[f'])\to\hom_{\pa\0}(\mr
f,\mr{f'})\to0.
\end{equation}
It follows from Lemma \ref{d=0} that
\begin{equation}\label{diz}
\hom_\ptr([f],[f'])=\hom_{\pa\0}(\mr f,\mr{f'})
\end{equation}
provided one of the following equations holds: $f=\id_X$, $f=!_X$,
$f'=\id_X$, $f'=^X!$, for an object $X\in\pa$. In particular
$$
\hom_\ptr(!_A,!_X)=\hom_{\pa\0}(!_A,!_X)=\hom_\pa(A,X).
$$
It follows that the full embedding $\pa\to\pa\0$ given by $X\mapsto\,!_X$ has
a unique lifting to $\ptr$.

\begin{Pro}\label{ras}
If $\pa$ is a pseudo-triangulated category, then for any object $X$ and for
any morphism $f:A\to B$ in $\pa$ one has the following exact sequences
$$
\cdots\to\hom_\pa(A[n+1],X)\to
\hom_\ptr(\tau^n([f]),!_X)\to\hom_\pa(B[n],X)\xto{f^*}\hom_\pa(A[n],X)\to\cdots
$$
and
$$
\cdots\to\hom_\pa(X,A[n-1])\xto{f_*}\hom_\pa(X,B[n-1])\to\hom_\ptr(^X!,\tau^n([f]))\to\hom_\pa(X,A[n])\to\cdots
$$
\end{Pro}

\begin{proof}
We prove exactness only for the first sequence. The proof for the second
sequence is similar and therefore we omit it. By the exact sequence \ref{des}
we have
$$
0\to \vd(\mr{f},!_X)\to \hom_\ptr([f],!_X)\to \hom_{\pa\0}(\mr{f},!_X)\to 0.
$$
It follows from the definition of the category $\pa\0$ that for
$f:A\to B$ one has the  exact sequence
$$
0\to\hom_{\pa\0}(\mr f,!_X)\to\hom_\pa(B,X)\xto{f^*}(A,X).
$$
This and Lemma \ref{d=0} imply exactness of the following sequence:
$$
\hom_\pa(B[1],X)\xto{f^*}\hom_\pa(A[1],X)\to
\hom_\ptr([f],!_X)\to\hom_\pa(B,X)\xto{f^*}\hom_\pa(A,X).
$$
Replacing $f$ by the translations of $f$ we get the result.
\end{proof}

\subsection{Homology}\label{sbs:homology}

We would like to introduce the notion of the homology in the setup
of pseudo-triangulated categories generalizing the classical notion
for triangulated categories. As in algebraic topology, a homology
must satisfy the exactness and excision axioms. To introduce these
axioms we need some preparations.

\begin{Le} Let $\pa$ be a pseudo-triangulated category. For any
objects $A$ and $B$ of $\pa$ one has a natural isomorphisms
$$\hom_\ptr(^A!, \,!_B)\cong \hom_\pa(A[1],B)$$
\end{Le}

\begin{proof} Since $\hom_{\pa^{[1]}}(^A!, \,!_B)=0$, the result
follows from the exact sequence (\ref{des}) and the fact that
$$\vu(^A!, \,!_B)=\vd(^A!, \,!_B)=\hom_\pa(A[1],B).$$

\end{proof}
In particular for any object $A$ there is a canonical morphism
$$j_A:^A!\to\,!_{A[1]}$$
corresponding to $\id_{A[1]}$. It follows from our construction that
$$p(j_A)=0.$$

Since $\vd(-,\,^A!)=0$, it follows that
$\hom_\ptr(-,\,^A!)=\hom_{\pa^{[1]}}(-,\,^A!)$. In particular for
any arrow $f:A\to B$ there is a canonical morphism
$$k_f:[f]\to \,^A!$$
in $\ptr$ corresponding to the morphism $(\id_A,0):\mr{f}\to ^A!$ in
$\pa^{[1]}$. By construction it is functorial in $[f]\in\ptr$.

 Let $f:X\to Y$ be a morphism in a
pseudo-triangulated category $\pa$. Then we have the following
commutative diagram in $\pa$
$$
\xymatrix{%
0\ar[r]\ar[d]&A\ar[d]^f\\
0\ar[d]\ar[r]&B\ar[d]^\id\\
A\ar[r]^f&B%
}
$$
which gives rise to the diagram in $\pa\0$
$$
!_A\xto{(0,f)}\,!_B\xto{(0,\id)}\,\mr f.
$$
Since $\vu(!_X,-)=0$, it has the unique lift
$$
!_A\xto{!_f}\,!_B\xto{i_f}\,[f]
$$
to $\ptr$. Let \begin{equation} \label{ravdakira} j_f:[f]\to
!_{A[1]}.
\end{equation} be the composite $j_f=j_A\circ k_f$.

Gluing these sequences and applying the translation functor we obtain the
sequence
$$
\cdots\to\,!_{A[n]}\to\,!_{B[n]}\to\,\tau^n([f])\to\,
!_{A[n+1]}\to\,\cdots
$$
which is functorial in $[f]\ptr$. Similarly one gets the sequence of
morphisms:
$$
\cdots\,\to{}^{A[n]}!\to\,^{B[n]}!\to\,\tau^{n+1}([f])\to\,
^{A[n+1]}!\to\, \cdots
$$

\begin{De}
A morphism $x:[f]\to[g]$ in $\ptr$ is called \emph{excising} if the induced
map
$$
\hom_\ptr(^X!,[f])\to \hom_\ptr(^X!,[g])
$$
is an isomorphism for any $X\in\pa$.
\end{De}

\begin{Le}\label{genUSA}
For any object $A$ the natural map
$$
j_A:\,^A!\to \,!_{A[1]}
$$
is excising.
\end{Le}

\begin{proof}
Since $\hom_{\pa\0}(^X!,\,!_Y)=0=\vu(^X!,\,^Y!)$ and
$\vu(^X!,\,!_Y)=\hom_\pa(X[1],Y)$ the result follows.
\end{proof}

Now we are ready to give the following definition.

\begin{De}
A \emph{homology} on a pseudo-triangulated category $\pa$ with values in an
abelian category $\mathsf A$ is a covariant functor $\h:\ptr\to\mathsf A$
satisfying the following two axioms:
\begin{itemize}
\item[(Exactness)]
For any morphism $f:A\to B$ of the category $\pa$ the sequence
$$
\h(!_A)\xto{!_f}\,\h(!_B)\xto{i_f}\,\h([f])\xto{j_f}\h(!_{A[1]})
$$
is exact.
\item[(Excision)]
If $x:[f]\to[g]$ is excising then $\h(x):\h([f])\to \h([g])$ is an
isomorphism.
\end{itemize}
\end{De}

In presence of the Excision Axiom, the Exactness Axiom is equivalent to the
assertion that for any $f:A\to B$ the sequence
$$
\h(^{B[-1]}!)\to \h([f])\to \h(^A!)\to \h(^B!)
$$
is exact. This easily follows from Lemma \ref{genUSA}.

For a homology $\h$ we put
$$
\h^n(A):=\h(!_{A[n]}), \ \ \h^n([f]):=\h(\tau^n([f])).
$$
Then we have an exact sequence
$$
\cdots\to\h^n(A)\to\h^n(B)\to\h^n([f])\to\h^{n+1}(A)\to\cdots
$$
natural in $[f]\in\ptr$.

\begin{Pro}
For any object $X\in\pa$ the functor
$$
\hom_\ptr(^X!,-):\ptr\to\Ab
$$
is a homology theory.
\end{Pro}

\begin{proof}
First we have to prove exactness of the
sequence
$$
\hom_\ptr(^X!,\,!_A)\to\,\hom_\ptr(^X!,\,!_B)\to\,
\hom_\ptr(^X!,\,[f])\to \hom_\ptr(^X!,\,!_{A[1]}).
$$
Since $\hom_{\pa\0}(^X!,\,Y_!)=0$ for all $Y\in \pa$, we have
$$
\hom_\ptr(^X!,\,!_Y)=\vu(^X!,\,!_Y)=\vd(^X!,\,!_Y)=\hom_\pa(X[1],Y),
$$
thanks to Lemma \ref{d=0} and Exact Sequence (\ref{des}). Now
exactness follows from Proposition \ref{ras}. It remains to prove
that the functor $\hom_\ptr( ^X!,-)$ transforms excising morphisms
to isomorphisms. But this is  obvious.
\end{proof}

\begin{Le}
Let $\ta$ be a triangulated category and $\e:\ta\to \Ab$ be a homology in the
classical sense. Then the functor $\h:\trn(\ta)\to\Ab$ defined by
$$
\h([f]):=\e(C_f)
$$
is a homology on the pseudo-triangulated category associated to the
triangulated category $\ta$. In this way one gets an equivalence
between the category of homologies in classical and new sense.
\end{Le}

\begin{proof}
By our definition of the category $\trn(\ta)$ the assignment
$f\mapsto C_f$ can be considered as a  well-defined functor
$\trn(\ta)\to\ta$. By Lemma \ref{conrep} a morphism $[f]\to[g]$ is
excisable iff the induced morphism $C_f\to C_g$ is an isomorphism.
From these facts, the first part of the statement follows.

Assume $\h$ is a homology in the new sense. For any morphism $f$ in $\ta$ the
morphism $(0,u_f):[f]\to \,!_{C_f}$ is excising. Hence $\h([f])=\e(C_f)$,
where $\e:\ta\to\Ab$ is given by $\e(A):=\h(!_A)$. It follows easily from
Exactness Axiom that $\e$ is a homology in the classical sense, hence the
result.
\end{proof}

\subsection{Massey triple product}

Let
$$
\xymatrix{
X\ar[r]^f&Y\ar[r]^g&Z\ar[r]^h&W
}
$$
be a diagram in a pseudo-triangulated category $\pa$. Suppose $hg=0$
and $gf=0$. Then we have the following commutative diagram in $\pa$:
$$
\xymatrix{
X\ar[r]\ar[d]_f&0\ar[d]^0\\
Y\ar[r]^g\ar[d]_0&Z\ar[d]^h\\
0\ar[r]&W
}
$$
which can be considered as the following diagram in $\pa\0$:
$$
^X!\xto{(f,0)} \mr{g}\xto{(0,g)}\,!_W.
$$
Observe that the composite morphism is zero in $\pa\0$. Since the functor
$p:\ptr\to \pa\0$ is identity on objects and surjective on morphisms the
diagram can be lifted to $\ptr$:
$$
[^X!]\xto x[g]\xto w\,[!_W],
$$
where $x$ and $w$ are morphisms of pseudo-triangles such that
$p(x)=(f,0)$ and $p(w)=(0,h)$. Then $p(wx)=0$, hence
$$
wx\in\vu(^X!,\,!_W);
$$
since $X\to0$ and $0\to W$ are split morphisms, the last groups can be
replaced by $\vd(^X!,\,!_W)$. Hence Lemma \ref{d=0} implies that
$$
wx\in \hom_\pa(X[1],W).
$$
Actually, this element depends on lifting. If one chooses $x_1$ and $w_1$
instead of $x$ and $w$, then we can write $x_1=x+a$ and $w_1=w+b$, where
$$
a\in\vu(^X!,\mr g)=\vd(^X!,\mr
g)=\cok(\hom_\pa(X[1],Y)\xto{g_*}\hom_\pa(X[1],Z)
$$
and
$$
b\in\vu(\mr g,!_W)=\vd(\mr
g,!_W)=\cok(\hom_\pa(Z[1],W)\xto{g^*}\hom_\pa(Y[1],W).
$$
It follows that $w_1x_1=wx+bx+wa$, therefore the class $\{h,g,f\}$ of $wx$ in
the quotient
$$
\frac{\hom_\pa(X[1],W)}{h_*\hom_\pa(X[1],Z)+f^*\hom_\pa(Y[1],W)}
$$
is invariant; we call it the \emph{Massey product}. By definition we
have $wx\in \{h,g,f\}$. In the case of triangulated categories it
coincides with the classical Massey product as defined in \cite{GM}.

The following fact is well-known  \cite[Theorem 13.2]{heller}.
\begin{Le}\label{toda1} Let
$$ \xymatrix{A\ar[r]^{a}& B\ar[r]^{b} &C\ar[r]^{c}
&A[1]}.
$$
be an acyclic triangle in a triangulated category. Then it is a
distinguished triangle if and only if
 $\id_{A[1]}\in \{c,b,a\}$.
\end{Le}

\subsection{$K_0$ for pseudo-triangulated categories} Let $\pa$ be  a
small pseudo-triangulated category. We let $K_0(\pa)$ be the abelian
group generated by the symbols $[X]$ where $X$ is an object of
$\pa$, modulo the relations K1-K3 below.
\begin{itemize}
\item[K1)] $[0]=0$,
\item[K2)] $[X]=[Y]$ provided there exists an isomorphism $f:X\to Y$ in
$\pa$,
\item[K3)] For any arrows $f:X\to Y$ and $f':X'\to Y'$ in $\pa$ and excising morphism
$x:[f]\to[f']$ in $\ptr$ one has $[X]+[Y']=[X']+[Y].$ \end{itemize}

One easily sees that this notion generalizes the Grothendieck's
original definition for triangulated categories.

\section{The class $\vartheta$ as the first obstruction}\label{inimage}

Recent work of Muro and his coauthors \cite{muro}, \cite{muroetc}
shows that not all triangulated categories have models. It turns out
that the class $\vartheta$ is the first obstruction for a
triangulated category to have a model. Namely we will prove that if
$\vartheta$ is not lies in the image of the canonical homomorphism
$\HH^2(\ta\0,\vd)\xto\th\HH^2(\ta\0,\vt)$ then $\ta$ has no models.
In other words we prove that if $\ta$ is a triangulated category
associated to a  stable model category or a Frobenious category then
 then the extension (\ref{hemigafa}) is a pushforward construction along the
transformation $\th:\vd\to\vt$ as it is defined in Section
\ref{domination} and hence the class $\vartheta\in\HH^2(\ta\0,\vt)$
lies in the image of the canonical homomorphism
$\HH^2(\ta\0,\vd)\xto\th\HH^2(\ta\0,\vt)$. Actually all this is an
easy consequence of the work of Baues \cite{TTC}, \cite{todabr}.

We also check that for the triangulated category constructed in
\cite{muro}, \cite{muroetc} the class $\vartheta$ does not lies in
the image of the homomorphism
$\HH^2(\ta\0,\vd)\xto\th\HH^2(\ta\0,\vt)$. This give an alternative
proof of the corresponding result of \cite{muro}, \cite{muroetc}.

\subsection{Push-forward construction and domination}\label{domination}

Let $\pa$ be a $\tau$-category equipped with a pseudo-triangulated
category structure given by a singular $\tau$-extension
$$
0\to\vu\xto i\ptr\xto p\pa\0\to0
$$
and a $\tau$-transformation $\vf:\vd\to\vu$.

Assume a $\tau$-transformation $\xi:\vu\to \vu_1$ of
$\tau$-bifunctors  is given which is an isomorphism as soon as one
of the arguments is a split morphism. Consider the following
category $\ptr_1$. The objects of $\ptr_1$ are the same as of the
categories $\pa\0$ and $\ptr$, i.~e. they are arrows of the category
$\pa$. Moreover $\hom_{\ptr_1}([f],[g])$ is defined using the
pushout  diagram of abelian groups:
$$
\xymatrix{
\vu(\mr f,\mr g)\ar[r]\ar[d]&\hom_\ptr([f],[g])\ar[d]\\
\vu_1(\mr f,\mr g)\ar[r]&\hom_{\ptr_1}([f],[g]) }
$$
It is  easy to see that in this way one gets a singular $\tau$-extension
structure on $\ptr_1$, such that the following diagram commutes:
$$
\xymatrix{
0\ar[r]&\vu\ar[r]\ar[d]^\xi&\ptr\ar[r]\ar[d]^j&\pa\0\ar[r]\ar[d]_\id&0\\
0\ar[r]&\vu_1\ar[r]&\ptr_1\ar[r]&\pa\0\ar[r]&0.
}
$$
Hence $\ptr_1$ together with the $\tau$-transformation
$\xi\circ\vf:\vd\to \vu_1$ is a pseudo-triangulated category
structure on $\pa$, called the \emph{pushforward construction}. In
this situation we also say that the pseudo-triangulated category
$\ptr$ dominates $\ptr_1$ and write $\ptr_1\le\ptr$.

The proof of the following easy fact is left to the reader.
\begin{Le} Massey triple product is invariant under dominations.
\end{Le}

\subsection{Toda triangulated categories} A  \emph{Toda triangulated category} is a triangulated category
$\ta$ such that the associated pseudo-triangulated category $\trn$
is dominated by a Toda pseudo-triangulated category. The following
is a straightforward.
\begin{Le} A triangulated category $\ta$ is a Toda triangulated
category iff the corresponding class $\vartheta\in \HH^2(\ta\0,\vt)$
lies in the image of the homomorphism
$\HH^2(\ta\0,\vd)\to\HH^2(\ta\0,\vt)$.
\end{Le}

\subsection{Track categories}\label{sec:TTC}

In the recent work \cite{TTC} Baues managed to construct Verdier
triangulated categories from the data which he called
\emph{triangulated track categories}. Recall that a \emph{track
category} $\B$ is a 2-category all of whose 2-morphisms are
invertible. Such categories appeared already in the classical work
\cite[Ch. V]{GZ}. Thus $\B$ consists of objects $X$, $Y$, etc., with
1-morphisms $\xi$, $\eta$ and with  2-morphisms $H:\xi\then \eta$.
If $\xi,\eta:X\to Y$ are 1-morphisms and there exists a 2-morphism
$H:\xi\then \eta$ then we say that $\xi$ and $\eta$ are homotopic.
The corresponding quotient category is denoted by $\B_\ho$, which
comes with the quotient functor $Q:\B\to \B_\ho$. Following
\cite{TTC} we use additive notation for the composite of
2-morphisms. A triangulated track category is a track category with
some extra data. We refer to the original paper of Baues \cite{TTC}
for the exact definition. Here we point out that any pointed
simplicial closed model category which is ``stable'' (meaning that
the suspension induces an auto-equivalence of the homotopy category)
gives rise to a triangulated track category structure on the track
category $\B$, which consists of fibrant-cofibrant objects,
1-morphisms are usual morphisms, while 2-morphisms are homotopy
classes of homotopies.


\subsection{Hardie category}

We need the following construction due to Hadrie \cite{hardie} which
we learned from \cite{todabr}. Let $\B$ be a  track category. Let
$\A$ be the corresponding homotopy category $\A=\B_\ho$. For each
morphism $f$ of $\A$ we choose its representative $\wt f$ in the
homotopy class of $f$. Hence $Q({\wt f})=f$.

Objects of the \emph{Hardie category} $\ha(\B)$ associated to the track
category $\B$ are  morphisms of $\A$. An object of $\ha(\B)$ corresponding to
a morphism $f$ is denoted by $\{f\}$. A morphism $\{f\}\to\{g\}$ in the
category $\ha(\B)$ corresponding to $f:A\to B$ and $g:X\to Y$ is an
equivalence class of triples $(\xi,\eta,H)$, where $\xi:A\to X$ and
$\eta:B\to Y$ are 1-morphisms of the track category $\B$, while $H$ is a
2-morphism $H:\eta\circ\wt f\then\wt g\circ\xi$. Two such triples
$(\xi,\eta,H)$ and $(\xi',\eta',H')$ are equivalent if there are 2-morphisms
$G:\eta'\then \eta$ and $K:\xi\then\xi'$ such that
$$
H'= \wt gK+H+\wt fG.
$$
Let $\{\xi,\eta,H\}$ be the equivalence class of $(\xi,\eta,H)$. Composition
in the Hardie category is given by
$$
\{\xi,\eta,H\}\circ\{\xi_1,\eta_1,H_1\}=\{\xi\xi_1,\eta\eta_1,\eta H_1+\xi
H\}.
$$

\begin{Le}\label{bpstr} Let $\B$ be a triangulated track category.
Then there is a well-defined functor $p:\ha(\B)\to \A\0$ which is
identity on objects and on morphisms is given by
$$
p\{\xi,\eta,H\}=(Q(\xi),Q(\eta))
$$
Moreover, if $\B$ is a triangulated track category, then $p$ is a part of a
singular $\tau$-extension
$$
0\to\vd\to\ha(\B)\xto p\A\0\to0.
$$
\end{Le}

\begin{proof}
This fact modulo notation is due to Baues \cite{todabr}. The
extension is the same as his linear extension \cite[Equation (2),
page 266]{todabr}, which is defined for much more general track
categories. The only thing to check is that for triangulated track
categories $D^\sharp$ in the notation of \cite{todabr} is the Toda
bifunctor. But this follows immediately from the definition of
$D^\sharp$ given in \cite[Equation (2.2)]{todabr} and the fact that
$D(X,Y)=\hom_\A(X[1],Y)$ for triangulated track categories, see
\cite[Equation (2.7)]{TTC}.
\end{proof}


\subsection{Pushforward construction in action}\label{apbaues}

Let $\B$ be a triangulated track category. By \cite{TTC} the
homotopy category $\A:=\B_\ho$ possesses a structure of triangulated
category and therefore we have a singular $\tau$-extension (see
Extension \ref{hemigafa}):
$$
0\to\vt\to\trn(\A)\to\A\0\to0.
$$
In this section we prove the following result.

\begin{Pro}\label{541tr}
Let $\B$ be a triangulated track category. Then there is a  functor
$T:\ha(\B)\to\trn(\A)$ which makes the diagram
$$
\xymatrix{
0\ar[r]&\vd\ar[r]\ar[d]_\th&\ha(\B)\ar[d]_T\ar[r]^p&\A\0\ar[r]\ar[d]_\id&0\\
0\ar[r]&\vt\ar[r]&\trn(\A)\ar[r]&\A\0\ar[r]&0
}
$$
commute. In particular for the class $\vartheta$ defined via Extension
(\ref{hemigafa}) one has
$$
\vartheta=\th_*(\beta),
$$
where $\beta\in\HH^2(\A\0,\vd)$ is the class of the extension constructed in
Lemma \ref{bpstr}.
\end{Pro}

\begin{proof} In the notations of \cite[Section 4]{TTC} the functor $T$ is defined by
$$T(\{f\})=(A\xto{f} B\xto{u} C_{{\wt f}}\xto{v} A[1])$$
where $f:A\to B$ is a morphism of $\A$ and  $u=Q(i_{\wt f})$ and $v=
Q(q_{\wt f})$.
\end{proof}

\subsection{Alternative approach}\label{apfranke} To obtain the previous result
that Extension (\ref{hemigafa}) for
a derived category of a differential algebra or a ring spectrum is
 pushforward along $\th$ instead of triangulated track categories we could have used systems of
 triangulated diagram categories in the sense of Franke \cite{franke}. In fact let  $\ka$ be a
 such system. In particular the categories $\ka_{\bf C}$ are given for any (finite) poset
 $\bf C$ satisfying some extra conditions. It follows from these axioms that each
 category $\ka_{\bf C}$ has a canonical structure of  a Verdier triangulated category.
These categories should be considered as refinement of the
triangulated category $\A=\ka_{\underline{0}}$, which is the base of
the system. Here   $\underline n$ denotes  the poset $\{0\leq \cdots
\leq n\}$. Based on the spectral sequence (32) \cite[Proposition
I.4.10]{franke} one can prove that there is a singular
$\tau$-extension
$$
0\to\vd\to\ka_{\underline{1}}
\to\A\to 0.
$$
and there is a  functor
$T:\ka_{\underline{1}}\to\trn(\A)$ which makes the diagram
$$
\xymatrix{
0\ar[r]&\vd\ar[r]\ar[d]_\th&\ka_{\underline{1}}
\ar[d]_T\ar[r]&\A
\0\ar[r]\ar[d]_\id&0\\
0\ar[r]&\vt\ar[r]&\trn(\A)\ar[r]&\A\0\ar[r]&0
}
$$
commute. The construction of the functor $T$ is similar to one
constructed in the proof of Proposition \ref{541tr} and is based on
the cones constructed in \cite[Section 1.4.6]{franke}.

It should be point out that if a triangulated category is associated
to a stable simplicial model category then both refinements --
triangulated track category  as well as system of triangulated
diagram categories are available. One can prove that in this case
Hardie category $\ha$ is equivalent to $\ka_{\underline{1}}$ and
hence both approach gives the same singular $\tau$-extensions.

\subsection{Muro's example}\label{muroexm}

For a small preadditive category $S$ we let ${\fa}(S)$ be the
additive completion of $S$. If $S$ has only one object (and hence
$S$ is just  a ring) then $\fa(S)$ is the category of finitely
generated free $S$-modules. Muro \cite{muro} shoved that the
category $\fa( \Z/4\Z)$ with the identity translation functor has
the unique triangulated category structure such that the triangle
$$\Z/4\Z\xto{2}\Z/4\Z\xto{2}\Z/4\Z\xto{2}\Z/4\Z$$
is distinguished. In this section we show that for this triangulated
category the extension \ref{hemigafa} is not a  pushforward along
$\th$. In the light of Section \ref{apbaues} and Section
\ref{apfranke}, it follows that this triangulated category does not
admits any refinement as a triangulated  track category \cite{TTC}
or as a system of triangulated diagram categories \cite{franke}.
This fact sharpers some results from \cite{muro}, \cite{muroetc}.

Consider the preadditive category $R$ which is generated by the
following graph
$$
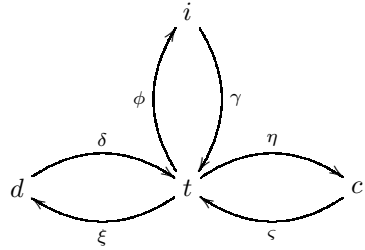
\xymatrix {&& i \ar@/^3ex/[dd]^{\gamma}&&\\
&&&&&\\
 d  \ar@/^3ex/[rr]^{\delta} && t \ar@/^3ex/[uu]^{\phi}
 \ar@/^3ex/[ll]^{\xi}\ar@/^3ex/[rr]^{\eta}&& c \ar@/^3ex/[ll]^{\varsigma}}
$$
modulo the following relations: all arrows are annihilated by $4$
and furthermore
$$2\cdot \delta \xi=0,\ \ \ 2\cdot \varsigma \eta=0,$$
$$\eta \delta=0, \ \phi\varsigma=0, \ \ \xi\gamma=0,$$
$$ \xi\delta=2\cdot \id_d, \ \ \eta\varsigma= 2\cdot \id_c, \ \ \phi \gamma=
2\cdot\id_i,$$
$$ \gamma \phi= \delta \xi+\varsigma \eta.$$
Then $\hom_{R}(d,c)=\hom_{R}(c,i)=\hom_{ R}(i,d)=\hom_{R}(d,c)=0$.
Moreover the abelian groups $\hom_ R(d,i)$, $\hom_{R}(d,t)$,
$\hom_{R}(c,d)$,  $\hom_{R}(c,t)$, $\hom_{R}(i,c)$, $\hom_{R}(i,i)$,
 $\hom_{R}(t,d)$, $\hom_{R}(t,c)$ and $\hom_{R}(t,c)$ are
isomorphic to $\Z/4\Z$. The rings $\hom_{R}(d,d)$, $\hom_{R}(c,c)$
and $\hom_{R}(i,t)$ are isomorphic to $\Z/4\Z$, while  $\hom_{
R}(t,t)$ as a ring is isomorphic to the ring
\begin{equation}\label{endt}
\hom_R(t,t)\cong \{(a,b,c)\in (\Z/4\Z)^3\mid a\equiv b \equiv c
\equiv 0{\rm (mod} 2)\}
\end{equation}
This isomorphism is given by
$$(2,2,0)\mapsto \gamma\phi, \ \ (0,2,2)\mapsto \delta\varsigma. $$

Let $R_1$ be the quotient of $R$ by the relations
$$2\xi=0, \ \  2\varsigma=0,\ \ \ \xi\varsigma=0.$$
Finally let $R_2$ be the quotient of $R_1$ by the relation
$$\gamma\phi=2\cdot \id_t.$$
We let $q:R\to R_2$ and $p:R_1\to R_2$ be the quotient
homomorphisms. We claim that neither $q$ and nor $p$ has a section.
This is clear for $p$ because even the homomorphism of abelian
groups $\Z/4\Z=\hom_{R}(t,s)\to \hom_{R_2}(t,s)=\Z/2\Z$ does not
have a section. For the functor $p$ one observes that
$p(x,y):\hom_{R_1}(x,y)\to \hom_{R_2}(x,y)$ is an isomorphism for
all possible $x,y$ except the case when $x=y=t$. Hence, if $p$ has a
section $s$, then $s$ would respects all arrows indicated in the
graph. But this contradicts to the fact that the equality
$\gamma\phi=2\cdot \id_t$ holds in $R_2$ but not in $R_1$.

Define $R_2$-$R_2$-bimodules $\vd,\vt,\vt_1$ as follows. The
bifunctor $\vt_1$ is zero everywhere but $\vt_1(t,t)=\Z/2\Z$. The
left and right action of the endomorphism ring of $t$ on
$\vt_1(t,t)$ is given by the multiplication on $a$ (which is the
same as the multiplication by $b$ or $c$). Here we used the
identification \ref{endt}. Moreover, we have
$$\vd(i,-)=0=\vd(-,i),\ \ \vd(d,-)=0=\vd(-,c),$$
$$\vd(c,d)=\Z/4\Z, \ \ \ \vd(t,d)=\vd(t,t)=\vd(c,t)=\Z/2\Z.$$
The arrows of $R_2$ acts on $\vd$ as follows. The homomorphisms
$\vd(c,\delta), \vd(t,\delta)$ are natural epimorphisms $\Z/4\Z\to
\Z/2\Z$, the morphisms $\vd(c,\xi), \vd(\varsigma,d)$ are natural
inclusions $\Z/2\Z\to \Z/4\Z$, finally we have
$\vd(t,\xi)=0=\vd(\varsigma,t)$, while $\vd(t,\delta)$,
$\vd(\eta,t)$ are isomorphisms. The bifunctor $\vt$ on objects has
the same values as the bifunctor $\vd$ and even morphisms act on
$\vt$ and $\vd$ in the same way provided the group $\vd(t,t)$ is not
involved. The rest actions are  given as follows. The morphisms
$\vt(t,\delta), \vt(t,\xi), \vt(\varsigma,t)$ are isomorphisms,
while $\vt(\eta,t)=0$. Then one has a binatural transformation
$\theta:\vd\to\vt$, such that $\theta(x,y)$ is the identity morphism
for all possible $x$ and $y$ except the case when $x=t=y$ and in
this exceptional case we have $\theta(t,t)=0$. One observes that we
have the following diagram with exact columns and rows
$$\xymatrix{&\vd\ar[d]_\theta &&&\\
0\ar[r]& \vt\ar[d] \ar[r]& R\ar[r]^{q}\ar[d]& R_2\ar[d]_{\id}\ar[r]  & 0\\
0\ar[r]& \vt_1 \ar[r]\ar[d]& R_1\ar[d]\ar[r]^{p}& R_2\ar[r]& 0\\
&0&0&&}$$ We have already seen that the bottom singular extension
does not split. Hence the middle singular extension is not a
pushforward along $\theta$.

All this related to Muro's example as follows. By mapping
$$d\mapsto (0\to \Z/4\Z),\ \ c\mapsto (\Z/4\Z\to 0), \ \
i\mapsto (\Z/4\Z\xto{1} \Z/4\Z),\ \ t\mapsto (\Z/4\Z\xto{2}
\Z/4\Z)$$
$$ \phi\mapsto (2,1), \eta\mapsto (1,0), \delta\mapsto (0,1),
\gamma\mapsto (1,2), \varsigma\mapsto (2,0), \xi\mapsto (0,2)$$ one
gets an equivalence of categories:
$$
\fa(\Z/4\Z)^{[1]}\cong \fa(R_2),$$
 while mapping
$$d\mapsto (0\to \Z/4\Z),\ \ c\mapsto (\Z/4\Z\to 0), \ \
i\mapsto (\Z/4\Z\xto{1} \Z/4\Z),\ \ t\mapsto (\Z/4\Z\xto{2}
\Z/4\Z)$$
$$ \phi\mapsto (2,1,0), \eta\mapsto (1,0,2), \delta\mapsto (0,1,2),
\gamma\mapsto (1,2,0), \varsigma\mapsto (2,0,1), \xi\mapsto
(0,2,1)$$ one gets an equivalence of categories: $$\trn\cong
\fa(R)$$ and these equivalences are compatible with bifunctors
$\vd,\vt,$
 etc. This proves that for the Muro's triangulated category the
extension \ref{hemigafa} is not a  pushforward along $\th$.

\section{Pseudo-triangulated versus triangulated categories}\label{muro}

\subsection{Embedding under domination}

Let us recall that we have a full embedding $!_?:\pa\to \pa\0$. Since
$\vu(!_X,-)=0$ this embedding has a unique lifting $!_?:\pa\to \ptr$ which is
still an embedding. Because of uniqueness it is invariant under domination.
In fact we have the following result.

\begin{Le}\label{le:51}
Let
$$
\xymatrix{
0\ar[r]&\vu\ar[r]\ar[d]^\xi&\ptr\ar[r]\ar[d]^j&\pa\0\ar[r]\ar[d]_\id&0\\
0\ar[r]&\vu_1\ar[r]&\ptr_1\ar[r]&\pa\0\ar[r]&0
}
$$
be part of a pushforward construction of pseudo-triangulated
categories. Then the diagram
$$
\xymatrix{
\pa\ar[r]^{!_?}\ar[dr]_{!_?'}&\ptr\ar[d]^j\\
&\ptr_1
}
$$
commutes. Moreover the functor $!_?:\pa\to \ptr$ has a left adjoint iff the
functor $!_?:\pa\to\ptr'$ does.
\end{Le}

\begin{proof}
The first part is a consequence of the uniqueness of lifting. To prove the
second part, we recall some general facts related to the adjoint functors.
Let $\c$ be a full subcategory of a category $\c_1$ and $x\in\c_1$. In these
circumstances one denotes by $x/\c$ the category of arrows $x\to c$, where
$c\in\c$. It is well known that the inclusion $\c\subset\c_1$ has a left
adjoint iff for all objects $x\in\c_1$ the category $x/\c$ has an initial
object.

According to Sequence \ref{des} we have the following commutative diagram
with exact rows
$$
\xymatrix{
0\to\vu(\mr{f},\,!_A)\ar[d]\ar[r]& \hom_\ptr([f],\,!_A)\ar[r]\ar[d]^j&\hom_{\pa\0}(\mr f,\,!_A)\ar[r]\ar[d]^\id&0\\
0\to\vu'(\mr f,\,!_A)\ar[r]&\hom_{\ptr'}([f],\,!_A)\ar[r]&\hom_{\pa\0}(\mr
f,\,!_A)\ar[r]&0.
}
$$
Since $\vu(-,!_A)\to\vu'(-,!_A)$ is an isomorphism, the middle vertical map
is also an isomorphism. It follows that for a fixed $f$ the category of
arrows $[f]\to !_A$ in $\ptr$ where $A$ runs over $\pa$, and the category of
arrows $[f]\to !_A$ in $\ptr'$, $A\in\pa$, are equivalent. From this, the
result follows.
\end{proof}

A similar fact is true for $^?!$ as well.

\subsection{The main result}

\begin{Th}\label{vaiteorema}
Let $\pa$ be a $\tau$-category equipped with a pseudo-triangulated
category structure given by a singular $\tau$-extension
$$
0\to\vu\xto i\ptr\xto p\pa\0\to0.
$$
Assume the functor  $!_?:\pa\to\ptr$ has a left adjoint functor
$L:\ptr\to \pa$ with counit of the adjunction $w_f:[f]\to
!_{L([f])}$, where $f:A\to B$  is a morphism in $\pa$. Declare a
triangle
$$
X\to Y\to Z\to X[1]
$$
to be distinguished provided there is a morphism $f:A\to B$ and a
commutative diagram in $\pa$
\begin{equation}\label{pre-tr1}
\xymatrix{A\ar[r]^f \ar[d]^a& B\ar[r]^{u_f}\ar[d]^b&
L([f])\ar[r]^{v_f}\ar[d]^c&A[1]\ar[d]^{a[1]}\\
X\ar[r]&Y\ar[r]&Z\ar[r]&X[1]}
\end{equation}
where $a,b,c$ are isomorphisms in $\pa$, $!_{u_f}=w_fi_f$, while
$v_f:L([f])\to A[1]$ is the unique morphism such that $j_f=v_f\circ
w_f$. Here $j_f$ and $i_f$ are maps in the sequence
$$
!_A\,\xto{!_f} \,!_B\,\xto{i_f} \, [f]\,\xto{j_f} !_{A[1]}
$$
constructed in Section \ref{sbs:homology}. With this class of
distinguished triangles  Axioms TR1,TR2,TR4,TR5 of triangulated
categories hold.

Moreover TR3 holds (and hence $\pa$ is a
triangulated category) iff the functor $[f]\mapsto L(f)[-1]$ is
right adjoint to the functor $\pa\to \ptr$ given by $X\mapsto ^X!$
and additionally
$$\{v_f,u_f,f\}=\id_{A[1]}$$ holds for all $f$.  If this is so then $\trn(\pa)\leq
\ptr$.
\end{Th}

\begin{proof}
The functor $L$ has the following universal property: for any
morphism $x:[f]\to !_X$ in $\ptr$ there exists a unique morphism
$g:L([f])\to X$ in $\pa$ such that $x=!_g\circ w_f$. Applying this
to  the sequence (referred as the \emph{pretriangle} corresponding
to $[f]$)
$$
!_A\,\xto{!_f} \,!_B\,\xto{i_f} \, [f]\,\xto{j_f} !_{A[1]}
$$
we see that $v_f$ indeed exists and is unique.
 Since $!_?$ is
full and faithful we have $L(!_A)=A$. By applying the functor $L$ to
the pretriangle we obtain the sequence
\begin{equation}\label{pre-tr1}
A\xto fB\xto{u_f}L([f])\xto{v_f}A[1],
\end{equation}
in $\pa$. Now we are in a position to check Axioms. The axioms TR1
and TR4 are obvious. Now we verify the axiom TR5. Consider a
commutative diagram
$$
\xymatrix{
A\ar[r]^f\ar[d]^a&B\ar[r]^{u_f}\ar[d]^b&L([f])\ar[r]^{v_f}&A[1]\ar[d]^{a[1]}\\
A'\ar[r]^{f'}&B'\ar[r]&L([f'])\ar[r]&A'[1].
}
$$
Then $(a,b):\mr f\to\mr{f'}$ is a morphism in $\pa\0$. Hence there exists a
morphism $x:[f]\to[g]$ in $\pa$ such that $p(x)=(a,b)$. By functoriality $x$
induces corresponding morphism on pretriangles
$$
\xymatrix{
!_A\ar[r]^{!_f}\ar[d]^{!_a}&!_B\ar[r]^{i_f}\ar[d]^{!_b}&[f]\ar[r]^{j_f}\ar[d]^x&!_{A[1]}\ar[d]^{!_{a[1]}}\\
!_{A'}\ar[r]^{!_{f'}}&!_{B'}\ar[r]&[f']\ar[r]&!_{A'[1]}
}
$$
Applying $L$ we finally get the commutative diagram
$$
\xymatrix{
A\ar[r]^f\ar[d]^a&B\ar[r]^{u_f}\ar[d]^b&L([f])\ar[r]^{v_f}\ar[d]^c&A[1]\ar[d]^{a[1]}\\
A'\ar[r]^{f'}&B'\ar[r]&L([f'])\ar[r]&A'[1]
}
$$
where $c=L(x)$.

To verify the axiom TR2, we take $f=\id_A$. By adjointness we have
$$
\hom_\ptr(\Id_A,\,!_B)=\hom_{\pa}(L(\Id_A),B).
$$
Since $\hom_{\pa\0}(\Id_A,!_B)=0$ and $\vu(\Id_A,-)=\vd(\Id_A,-)=0$, it
follows that $\hom_\ptr(\Id_A,!_B)=0$ for all $B\in\pa$. Now the Yoneda lemma
shows that $L(\Id_A)=0$. It follows that
$$
A\xto{\id}A\to 0\to A[1]
$$
is a distinguished triangle.

If Axiom TR3 holds, then $\{v_f,u_f,f\}=\id_{A[1]}$ by Lemma
\ref{toda1}. It follows  from Lemma \ref{conrep} and analogue of
Lemma \ref{le:51} for right adjoint functors that the functor
$[f]\to L(f)[-1]$ is indeed the right adjoint functor.

Conversely, assume $\{v_f,u_f,f\}=\id_{A[1]}$  holds for all $f$ and
the the functor $[f]\to L(f)[-1]$ is  the right adjoint to the
functor $?^!$. We have to check Axiom TR3.

We start with observation that the functor $L$ takes any excising
morphism into isomorphism. In fact if $x:[f] \to [g]$ is a excising,
then all morphisms
$$\hom_\pa(X,L(f)[-1])\to \hom_\ptr(^X!,[f])\to \hom_\ptr(^X!,[g])\to \hom_\pa(X,L(g)[-1])$$
are isomorphism for all object $X\in \ta$. It follows from the
Yoneda lemma that $L(f)\to L(g)$ is also an isomorphism.

Next, remark that for any $f:A\to B$ there are morphisms $x:^A!\to
[u_f]$ and $w:[u_f]\to !_{A[1]}$ in $\ptr$ such that $p(x)= (f,0)$,
$p(w)=(0,v_f)$ and $wx\in \Upsilon(^A!,
!_{A[1]})=\hom_\pa(A[1],A[1])$ represents the identity morphism
$\id_{A[1]}$. Hence $wx=j_A$ is an excising morphism. It follows
that $L(wx)$ is an isomorphism. Thus $L(w):L([u_f])\to A[1]$ is a
split epimorphism. By 5-Lemma applied to exact sequences induced by
$u_f$ and $f$ it follows that  $L(w)$ is in fact an isomorphism.
This fact implies TR3. Hence $\pa$ is a triangulated category. The
proof also shows that the triangulated category structure is
dominated by $\ptr$.
\end{proof}
Now we are in a position to prove our main result.
\begin{Co}
Let $\pa$ be a $\tau$-category equipped with a pseudo-triangulated
category structure given by a singular $\tau$-extension
$$
0\to\vu\xto i\ptr\xto p\pa\0\to0.
$$
Then the following conditions are equivalent
\begin{itemize}
\item[i)]
There is a triangulated category structure $\trl(\pa)$ on $\pa$ and
a domination
$$
\trn(\pa)\le\ptr.
$$
\item[ii)]
There is a functor $L:\ptr\to \pa$ which is left adjoint to the
functor $X\mapsto !_X$, while $[-1]\circ L$ is a right adjoint to
the functor $X\mapsto ^X!$ and $\{v_f,u_f,f\}=\id_{A[1]}$ for all
$f$.
\end{itemize}
\end{Co}

\begin{proof} The implication i) $\then$ ii) follows from  Lemma \ref{conrep} the functor
together with Lemma \ref{le:51} and Lemma \ref{toda1}. The
implication ii) $\then$ i) follows from Theorem \label{vaiteorema}.
\end{proof}

\section{Idempotent completion}\label{bsch}
\subsection{ Karoubization}

Any additive category $\A$ has a Karoubian completion $\A\kar$,
which is a Karoubian category with a full embedding $i:\A\to \A\kar$
satisfying the following property. If $\B$ is a Karoubian category
and $j:\A\to \B$ is an additive functor, then there exists an
essentially unique functor $f:\A\kar\to\B$ with $j=fi$. Objects of
$\A\kar$ are pairs $(A,e)$, where $A$ is an object of $\A$ and
$e:A\to A$ is an idempotent. A morphism $(A,e)\to (A',e')$ is a
morphism $f:A\to A$ in $\A$ such that
$$
fe=e'f=f.
$$
Let us observe that the identity morphism of $(A,e)$ is $e$ and the
functor $i$ is given by $i(A)=(A,\id_A)$.

\begin{Le}\label{spltidcr} An idempotent $e$ of the category $\A$ is split iff $(A,e)$ as an object of
$\A\kar$ is isomorphic to an object of the image of the functor
$i:\A\to\A\kar$.
\end{Le}

\begin{proof} One easily checks that having mutually inverse morphisms $$a:(A,e)\to
(B,\id_B), \ \ \ b:(B,\id_B)\to (A,e)$$ is exactly the same as to
have a splitting data for $e$.
\end{proof}

\subsection{Lifting of idempotents}
It is well know that if $R\to S$ is a surjective homomorphism of
rings with nilpotent kernel then any idempotent of $S$ is an image
of an idempotent of $R$. We can specialize this for the ring
homomorphism $\hom_\A(A,A)\to \hom_{\A/\I}(A,A)$ to get the
following fact.
\begin{Le}\label{liftid}
Let $\I$ be a nilpotent ideal of an additive category $\A$. For any
idempotent $f:A\to A$ in the quotient category $\A/\I$ there is an
idempotent $e:A\to A$ in $\A$ such that $Q(e)=f$.
\end{Le}

\subsection{Singular extensions and idempotent completion}
Let
$$
0\to D\xto i\A\xto F\B\to 0
$$
be a singular extension of an additive category $\A$ by a biadditive
bifunctor $D:\A\op\x\A\to \Ab$. In this section we compare the
categories $\A\kar$ and $\B\kar$. Actually our results here are very
particular case of much more general results obtained in
\cite{triest}.

To make notations simpler we write $\bar{f}:A\to B$ instead of
$F(f):A\to B$. Here  $f:A\to B$ is amorphism in $\A$. Since all
idempotents in $\Ab$ splits the bifunctor $D$ has the canonical
extension $D\kar:(\A\kar)\op\x\A\kar\to \Ab$ to the category
$\A\kar$. In more details
$$D\kar((A,\bar{e}),(A',\bar{e}'))=\im(\bar{e}'_*\circ\bar{e}^*:D(A,A')\to D(A,A')).$$
This works because any idempotent in $\B$ has the form $\bar{e}$ for
an idempotent $e$ in $\A$ thanks to Lemma \ref{liftid}. Now we fix
idempotents $e:A\to A$, $e':A'\to A'$.  Since
$$\hom_{\A\kar}((A,e),(A',e'))=\im(e'_*\circ e^*:\hom_\A(A,A')\to \hom_\A(A,A'))$$
the inclusion $i:D(A,A')\to \hom_\A(A,A')$ has the unique extension
$$i\kar:D\kar((A,\bar{e}),(A',\bar{e}'))\to
\hom_{\A\kar}((A,e),(A',e')).$$ This allows to consider $D\kar$ as a
square zero  ideal in $\A\kar$. The corresponding quotient category
is denoted by $\tilde{\B}$. The  objects of the category
$\tilde{\B}$ are pairs $(A,e)$, where $e$ is an idempotent in the
category $\A$. The morphisms are
$$\hom_{\tilde{{B}}}((A,e),(A',e'))=\im(\bar{e}'_*\circ\bar{e}^*:\hom_\B(A,A')\to
\hom_\B(A,A')).$$ This follows from the fact that  the composites
$e'_*\circ e^*$ and $\bar{e}'_*\circ\bar{e}^*$ are idempotents and
hence  their images are in fact direct summands. It follows that the
functor $\tilde{B}\to \B\kar$ defined on objects by $(A,e)\mapsto
(A,\bar{e})$ is full and faithful and in fact an equivalence thanks
to Lemma \ref{liftid}. Having this equivalence in mind the bifunctor
$D\kar$ can be considered as a bifunctor on $\tilde{B}$. We can now
summarize our discussion.

\begin{Le}\label{triestagaf} If
$$
0\to D\xto i\A\xto F\B\to 0
$$
is a singular extensions of additive categories, then we have a
singular extension of additive categories
$$
0\to D\kar\xto i\A\kar\xto{\tilde{F}}\tilde{\B}\to 0
$$
and an equivalence of categories $\tilde{\B}\to \B\kar$. Moreover,
we have also a singular extension of additive categories
$$0\to D\kar\to \tilde{A}\to \B\kar\to 0$$
and an equivalence of categories
$$\tilde{A}\to \A\kar.$$

\end{Le}

\begin{proof} Only the second part of the statements needs some
comments. It follows from the first part  by notice that an
equivalence of categories yields an isomorphism in the
Baues-Wirsching cohomology \cite{BW}.

\end{proof}

Lemma \ref{triestagaf} says that up to equivalence of categories a
singular extension gives rise to a singular extension by passing
trough the idempotent completion. Based on this fact we now prove
the following easy fact.

\begin{Pro}\label{lechen} Let $$
0\to D\xto i\A\xto F\B\to 0
$$
be a singular extensions of additive categories and let $e:A\to A$
be an idempotent in $\A$. Then $e$ splits   iff $F(e)$ splits.

\end{Pro}

\begin{proof} Of course any functor takes split idempotents to split
ones. Assume now $\bar{e}=F(e)$ is split. This means that  there is
an isomorphism $(A,\bar{e})\to (B,\id_B)$ in $\B\kar$ (see Lemma
\ref{spltidcr}). But both  $(A,\bar{e})$ and $(B,\id_B)$ are in the
image of the functor $\tilde{\B}\to \B\kar$, which is an equivalence
of categories (see Lemma \ref{triestagaf}).  Hence there is an
isomorphism $\tilde{x}:(A,e)\to (B,\id)$ in $\tilde{\B}$. The
functor $\tilde{F}:\A\kar\to \tilde{\B}$ is full and reflects
isomorphisms thanks to Lemma \ref{triestagaf} and Lemma
\ref{refisonil}. It follows that there is an isomorphism $x:(A,e)\to
(B,\id_B)$ and the result follows.
\end{proof}
\bigskip

\subsection{Split idempotents in the category of arrows} Let $f:A\to B$ be a morphism of
an additive category $\A$. It is clear that a morphism
$(a,b):\mr{f}\to \mr{f}$ in the category $\A\0$ is an idempotent iff
$a:A\to A$ and $b:B\to B$ are idempotents.
\begin{Le}\label{idisari} An idempotent $(a,b):\mr{f}\to \mr{f}$  of $\A\0$ splits
iff $a$ and $b$ are split idempotents of the category $\A$.
\end{Le}
\begin{proof} Assume $a$ and $b$ are split idempotents of the category
$\A$. Let $A\xto{c} C\xto{d} A$ and $B\xto{s} D\xto{t} B$ be
splitting data for $a$ and $b$. We set $g=sfd:C\to D$. One easily
checks that $\mr{f}\xto{(c,d)} \mr{g}\xto{(d,t)} \mr{f}$ is a
splitting data of the idempotent $(a,b):\mr{f}\to \mr{g}$. The
converse statement is obvious.
\end{proof}

\subsection{Split idempotents in the category of pseudo-triangles}
Let $\pa$ be a $\tau$-category equipped with a pseudo-triangulated
category structure given by a singular $\tau$-extension
$$
0\to\vu\xto i\ptr\xto p\pa\0\to 0.
$$
and $\tau$-transformation $\varphi:\vd\to \vu$.
\begin{Le} Let $x:[f]\to [f]$ be an idempotent in $\ptr$ with
$p(x)=(a,b):\mr{f}\to \mr{f}$. Then $x$ is a split idempotent  in
$\ptr$ iff $a$ and $b$ are split idempotents in $\pa$.
\end{Le}

\begin{proof} If part is clear. Assume $a$ and $b$ are split
idemotents. Then $(a,b):\mr{f}\to \mr{f}$ is also split idempotent
thanks to Lemma \ref{idisari}. Hence the result follows from
Proposition \ref{lechen}.
\end{proof}

As an immediate consequence of the above abstract non-sense we get
the following crucial lemma in \cite{LeChen}.

\begin{Co}\label{chinid} Let $\ta$ be a triangulated category and let $$
\xymatrix{
A\ar[r]^{f}\ar[d]^{a}& B\ar[r]^{u_f}\ar[d]^{b} &C_f\ar[r]^{v_f}\ar[d]^{c} &A[1]\ar[d]^{a[1]}\\
A\ar[r]^{f}& B\ar[r]^{u_{f}} &C_{f}\ar[r]^{v_{f}} &A[1] }
$$
be a morphism of distinguished triangles. Assume $a,b$ and $c$ are
idempotents. If $a$ and $b$ are split idempotents, then
$(a,b,c):[f]\to [g]$ is a split idempotent of the category $\trn$.
In particular $c$ is a split idempotent of $\ta$.
\end{Co}
\subsection{The full embedding $\varrho:  ({\A\kar})\0\to (\A\0) \kar .$}
In this section we compare  categories $({\A\kar})\0$ and $(\A\0)
\kar $.

Let $\A$ be an additive category. Objects of $({\A\kar})\0$ are
arrows $f:(A,e)\to (A',e')$ in $\A\kar$, where $A$ and $A'$ are
objects of $\A$, while $e$ and $e'$ are idempotents of the category
$\A$. We can also say that the objects of the category
$({\A\kar})\0$ are  diagrams in $\A$
$$\xymatrix{
A\ar[r]^e\ar[d]_f\ar[dr]^f& A\ar[d]^f\\
A'\ar[r]_{e'}&A'}$$ such that $e^2=e,(e')^2=e',fe=f=e'f$. Such an
object is denoted by $(A,e,A',e',f)$.

On  the other hand the objects of $(\A\0) \kar $ are pairs
$(\mr{f},x)$, where $f:A\to B$ is an arrow in $\A$ and
$x=(e,e'):\mr{f}\to \mr{f}$ is an idempotent in $\A\0$. We can also
say that the objects of the category $(\A\0) \kar $ are diagrams in
$\A$
$$\xymatrix{
A\ar[r]^f\ar[d]_e& A'\ar[d]^{e'}\\
A\ar[r]_{f}&A'}$$ with $e^2=e,(e')^2=e',fe=e'f$. Such an object is
denoted by $[A,f,A',e,e']$. We see that the map
$$\varrho((A,e,A',e',f))=[A,f,A',e,e']$$
yields an embedding of the class of objects of $({\A\kar})\0$ into
the class of objects of $(\A\0) \kar $. The next lemma shows that
$\varrho$ can be extended as a full and faithful functor
$$\varrho:  ({\A\kar})\0\to (\A\0) \kar .$$

\begin{Le} For any objects $(A,e,A',e',f)$ and $(B,d,B',d',g)$ of
the category $ ({\A\kar})\0$ we have
$$\hom_{ ({\A\kar})\0}((A,e,A',e',f), (B,d,B',d',g))\cong \hom_{(\A\0) \kar
}([A,f,A',e,e'], [B,g,B', d,d'])$$
\end{Le}

\begin{proof} A direct inspection shows that in both cases morphisms
are   pairs $(h,h')$, where $h:A\to B$ and $h':A'\to B'$ are
morphisms in $\A$ such that
$$dh=h=hc, \ \ \ d'h'=h'=h'c', \ \ \ h'f=gh.$$
\end{proof}
\subsection{Idempotent completion of pseudo-triangulated categories}
In this section we show that Karubization of  a pseudo-triangulated
category carries a natural pseudo-triangulated category structure
(compare with \cite{ictc}). This is based on the previous
relationship between the categories $({\A\kar})\0$ and $(\A\0) \kar
$.

 Let $\pa$ be a
$\tau$-category equipped with a pseudo-triangulated category
structure given by a singular $\tau$-extension
$$
0\to\vu\xto i\ptr\xto p\pa\0\to 0.
$$
and $\tau$-transformation $\varphi:\vd\to \vu$. By passing to the
idempotent completion we obtain another singular $\tau$-extension
(see lemma \ref{triestagaf}):
$$
0\to{\vu}\kar\to\tilde{\ptr}\to{(\pa\0)}\kar\to  0.
$$
Now we can pull-back it along the $\varrho$ to get a singular
$\tau$-extension
$$
0\to{\vu}\kar\to\widehat{\ptr}\to{(\pa\kar)}\0\to 0
$$
which in fact is a pseudo-triangulated category structure on
$\pa\kar$. One easily sees that for triangulated categories this is
exactly the  construction in \cite{ictc}.

\

\

\noindent{\bf Acknowledgments}

\noindent Dieser Artikel wurde zwar nicht in Bonn geschrieben, ist
aber sehr vom mathematischen Umfeld in Bonn gepr\"agt, von dem ich
zu verschiedenen Stadien meiner mathematischen Laufbahn profitieren
konnte. Besonderen Einfluss hatten die Artikel von Hans Joachim
Baues, Jens Franke und Stefan Schwede. Sehr hilfreich waren
ausserdem die Vortr\"age von Mamuka Jibladze und Fernando Muro in
Bonn.


\begin{thebibliography}{33}

\bibitem{ictc}
{\sc P. Balmer} and {\sc M. Schlichting}. Idempotent completion of
triangulated categories. J. of Algebra, 236(2001), 819--834.

\bibitem{TTC}
{\sc H.-J. Baues}. Triangulated track categories. Georgian Math. J. 13 (2006), 607--634.

\bibitem{todabr}
{\sc H.-J. Baues}. On the cohomology of categories, universal Toda
brackets and homotopy pairs. $K$-Theory 11 (1997), no. 3, 259--285.

\bibitem{shukla}
{ \sc H.-J. Baues} and {\sc T. Pirashvili}. Comparison of Mac~Lane, Shukla and
Hochschild cohomologies. \emph{J. Reine und Angew. Math.}  598 (2006), 25--69.

\bibitem{BW}
{\sc H.-J. Baues} and {\sc G. Wirsching}. Cohomology of small
categories. \emph{J. Pure Appl. Algebra} 38 (1985), 187--211.

\bibitem{franke}
{\sc J. Franke}. Uniqueness theorem for certain triangulated
categories possessing Adams spectral sequence. Preprint, Bonn, 1996.
Avaialable at {\bf http://www.math.uiuc.edu/K-theory/0139/}.


\bibitem{GZ}
{\sc P. Gabriel} and {M. Zisman}. Calculus of fractions and homotopy theory.
Ergebnisse der Mathematik und ihrer Grenzgebiete, Band 35 Springer-Verlag New
York, Inc., New York 1967 x+168 pp.

\bibitem{GM}
{\sc S.I. Gelfand} and {\sc Y.I. Manin}. Methods of homological
algebra. Second edition. Springer Monographs in Mathematics.
Springer-Verlag, Berlin, 2003. xx+372 pp.

\bibitem{hardie}
{\sc K. A. Hardie}. On the category of homotopy pairs. Topology Appl. 14(1982), 59--69.


\bibitem{heller} A. Heller. Stable homotopy categories. Bull. Amer. Math. Soc. 74(1968),
28-63.

\bibitem{JP}
{\sc M. Jibladze} and {\sc T. Pirashvili.} Cohomology of algebraic theories.
J. of Algebra. 137(1991), 253-296.


\bibitem{LeChen}

{\sc J. Le} and {\sc X.W. Chen}. Karoubianness of a triangulated
category. J. of Algebra. 310(2007), 452-457.

\bibitem{Mit}
{\sc B. Mitchell}.  Rings with several objects. Advances in Math. 8,
(1972), 1--161.

\bibitem{muro}
{\sc F. Muro}. A triangulated category without models. {\bf arXiv:
math.KT/0703311}.

\bibitem{muroetc}
{\sc F. Muro}, {\sc S. Schwede} and {\sc N. Strickland}.
Triangulated category without models. {\bf arXiv:0704.1378}.

\bibitem{triest}
{\sc T. Pirashvili}. Projectives are free for nilpotent algebraic
theories.
 Algebraic $K$-theory and its applications (Trieste, 1997), 589--599,
 World Sci. Publ., River Edge, NJ, 1999.


\bibitem{puppe}
{\sc D. Puppe}. On the structure of stable homotopy theory. Colloquium on
algebraic topology. Aarhus Universitet Matematisk Institut (1962), 65--71.

\bibitem{sh}
{\sc U. Shukla}. Cohomologie des alg\'ebres associatives. \emph{Ann.
Sci. \'Ecole Norm. Sup.} (3) \textbf{78} (1961), 163--209.

\bibitem{verdier}
{\sc J. L. Verdier.} Des Cat\'egories deriv\'ees des cat\'egories
ab\'eliennes. Ast\'erisque  v. 239 (1996), 253 pp.

\end{thebibliography}
\end{document}